\documentclass[10pt,amssymb]{article}
\usepackage{amsfonts,latexsym}
 
\newcommand{\R}{\mathbb R}

\newcommand{\C}{\mathbb C} 
\newcommand{\del}{\partial}
\newcommand{\e}{\varepsilon}

\def\ds{\displaystyle}

\newtheorem{theorem}{Theorem} 
\newtheorem{lemma}{Lemma} 
\newtheorem{proposition}{Proposition} 
 
\newtheorem{definition}{Definition} 
\newtheorem{remark}{Remark} 

\begin{document} 
 
\title{Complete embedded minimal $n$-submanifolds in ${\C}^{n}$}
 
\author{Claudio Arezzo and Frank Pacard}
\date{}
\maketitle 

\section{Introduction}

A classical problem in the theory of minimal submanifolds of Euclidean 
spaces is to study the existence of a minimal submanifold with a prescribed behavior at
infinity, or to determine from the asymptotes the geometry of the whole submanifold. Beyond
the intrinsic interest of these questions, they are also of crucial importance when studying
the possible singularities of minimal submanifolds in general  Riemannian manifold. When
studying minimal surfaces, i.e. two  dimensional submanifolds, the standard tool to solve
these problems is given by the Weierstrass representation formula which  relates the geometry
of the minimal surface to complex analytic  properties of holomorphic one-forms on Riemann
surfaces. Recently, gluing technics have been developed and have provided an abundant number
of new examples of minimal hypersurfaces in Euclidean space.

\medskip

For higher dimensional minimal submanifolds such a link clearly  
disappears and no complex analysis can be put into play. While, gluing technics have been
extensively used in the study of  minimal hypersurfaces, they have not been adapted to handle
higher codimensional submanifolds. The aim of this paper is to use a gluing technique for
minimal submanifolds to make a  step towards the understanding of these questions in
arbitrary codimension. We will restrict ourselves to the  case of real $n$-dimensional
submanifolds of ${\C}^{n}$. There are two main  reasons for doing so. The first one is
technical: when  trying to desingularize the intersection, for  example, of a pair  of
$n$-planes, which give the desired asymptotic behavior,  one needs a model of a minimal
submanifold with this behavior at  infinity, to rescale and to glue into the pair of planes
where a neighborhood of the intersection is removed. This local model  needs to be
sufficiently simple to allow a very detailed study of the  linearized mean curvature
operator. In our situation this is provided by a generalization of an area minimizing
submanifold  found by Lawlor in \cite{Law-G}, while in more general cases such an example is
not known.  The second reason is that among minimal $n$-submanifolds of  ${\C}^n$ there is a
special family, namely the  special Lagrangian ones, which are of great importance in a
variety of geometric and physical problems (see, for  example,  \cite{Har-Law}, \cite{jo}
and  \cite{syz}).

\medskip

To better describe our result let us first observe that for  minimal surfaces in ${\C}^2$ 
graphs of meromorphic  (or anti-meromorphic) functions  are enough to answer some the above
questions. 

\medskip

For example, given $z_1, \ldots, z_{k} \in {\C}$   
 all distinct and $\alpha_1, \ldots,  \alpha_k \in {\C}$ , 
the surface which is the graph of  
\begin{equation}
\label{eq:minc2}
z \in {\C}\setminus \{z_1, \ldots, z_k\}\longrightarrow \alpha_1   
\,  (\bar{z} -\bar{z}_1)^{-1}+\ldots  + \alpha_k \, (\bar{z}-\bar{z}_k)^{-1} \in {\C}^2
\end{equation}

is a complete embedded minimal surface with $k+1$ ends.

\medskip

When $k=1$ we get the usual hyperbola 

$$z \in {\C}\setminus \{0\} \longrightarrow (z , \bar{z}^{-1}) \in  {\C}^2.$$

which can be seen as the two dimensional element of the following  family of  $n$-submanifolds, $H_I$, of  ${\C}^{n}$,

$$\left( 0, \frac{\pi}{n} \right) \times  S^{n-1} \ni (s,  \Theta) \longrightarrow 
\frac{e^{i  s}}{(\sin (ns))^{\frac{1}{n}}} \, \Theta \in  {\C}^{n}$$

For any $n$, $H_I$ is a complete special Lagrangian  (and therefore minimal) submanifold of  
${\C}^{n}$  with two ends \cite{Har-Law},  \cite{Law-B} and \cite{Law-G}.

\medskip

\noindent

Observe that, for any ${\cal R} \in O(n)$, we can define  $H_{\cal R}$ the 
$n$-submanifold of  ${\C}^{n}$ which is parameterized by

$$\left( 0, \frac{\pi}{n} \right) \times  S^{n-1} \ni (s, \Theta)  \longrightarrow  
\frac{1}{(\sin (ns))^{\frac{1}{n}}} \, (\cos s \,  \Theta + i \,  \sin s \, {\cal R} 
(\Theta)
) \, \in{\C}^{n}.$$

This is still a minimal submanifold of ${\C}^{n}$, since  
$
\left( \begin{array}{rllll}
I & 0\\
0 & {\cal R}
\end{array}
\right)
$
belongs to $O(2n)$. However, given the complex structure
in ${\R}^{2n}$ induced by  the  identification of
${\R}^{2n}$ with ${\C}^n$, it fails to be Lagrangian,
except when
${\cal R}^2 = I$.

\medskip

In this paper we prove the existence  of  minimal embedded  $n$-submanifolds 
having $k+1$ ends which, in some sense are  the high dimensional analogues of the  surfaces 
described by (\ref{eq:minc2}). Consider the  "horizontal" $n$-plane 

$$\Pi_0 : = \left\{ x \in {\C}^{n}\, : \, x\in {\R}^n \right\}$$

and, for all $j=1, \ldots , k$, consider the $n$-planes

$$\Pi_j : = \left\{ \left( x_j+ \cos \frac{\pi}{n}\, x + i \, \sin \frac{\pi}{n}\,  \, 
{\cal R}_j \,    x \right) \in {\C}^{n}\, :  \, x \in {\R}^n \right\},$$

where $x_1, \ldots, x_{k}\in {\R}^n$ and ${\cal R}_1, \ldots,  {\cal R}_{k} \in O(n)$ 
are fixed.  We further assume that all $x_j$ are distinct. One property of  this set of
planes is that, than all  angles between the planes $\Pi_0$ and $\Pi_j$, are equal to 
$\frac{\pi}{n}$. We refer to \cite{Har-Law} for a definition of the angles between pairs of
$n$-planes  in ${\mathbb C}^n$.

\medskip

We set

\[
\xi_{j,j'}: = \frac{x_j-x_{j'}}{|x_j-x_{j'}|},
\]

if $j\neq j'$ and we define the $k \times k$ matrix  $\Gamma := (\gamma_{jj'})_{j,j'}$ 
by $\gamma_{jj}=0$  for all $j=1, \ldots, k$ and

\[
\gamma_{jj'} :=  \frac{1}{|x_j-x_{j'}|^{n}} \, {\int}_{S^{n-1}}  \left( {\cal R}_j  
\Theta \cdot {\cal R}_{j'} \Theta  - n \, (\Theta \cdot {\cal R}_j  \xi_{j,j'})\, (\Theta
\cdot {\cal  R}_{j'} \xi_{j,j'}) \right) \, d\theta ,
\]

for all $j \neq j'$.

\medskip

Finally, for any ${\cal A}_0 \in M_n ({\R})$, we define for all $j=1, \ldots , k$

\[
\lambda_j : = -  \int_{S^{n-1}} {\cal A}_0 \Theta \cdot
{\cal R}_j  \Theta \, d \theta ,
\]

and $\Lambda : = (\lambda_1, \ldots, \lambda_{k})$.

\medskip

In order for the desingularization to be possible, we will do the following assumptions :

\medskip

\begin{itemize}

\item[(H1)] $ \qquad \qquad \qquad \qquad  \forall j \neq j' , \qquad  \xi_{jj'}\notin 
\mbox{Im} \,  (I - {\cal R}_{j'}^{-1}\, {\cal R}_j) $

\item[(H2)] $\qquad \qquad \qquad \qquad \mbox{The matrix  $\Gamma$ is invertible}$

\item[(H3)] $\qquad \qquad \qquad \qquad  \Gamma^{-1}\, \Lambda  \in (0,+\infty)^{k}$

\end{itemize}

\medskip

\noindent

These assumptions will be commented in the next section.  We can now state the main 
result of the paper

\begin{theorem}

Assume that $n\geq 3$, and $(H1)$, $(H2)$ and $(H3)$ are fulfilled. Then, the set  of 
$n$-planes  
$\Pi_0, \ldots, \Pi_{k}$ can be desingularized to produce a  complete,  embedded  minimal
$n$-dimensional  submanifolds of ${\C}^{n}$ which has $k+1$ planar ends. This minimal   
$n$-submanifold is of the topological type of a $n$-sphere with $k+1$   punctures and has
finite  total curvature.

\medskip

More precisely, there exists $\e_0$ and for all $\e \in (0, \e_0)$  there exists 
$\Sigma_\e$ a  $n$-dimensional submanifold of ${\C}^{n}$ such that 

\begin{enumerate}

\item[(i)] For all $\e\in (0, \e_0)$,  $\Sigma_\e$ is embedded and minimal.

\item[(ii)] For all $\e\in (0, \e_0)$, $\Sigma_\e$ has $k+1$ ends,  which, up to
 translations,  are given by $\Pi_j$,  for $j=1, \ldots, k$ and $\{ x + i \, \e  \,  {\cal
A}_0 \, x \, | \, x \in {\R}^n\}$.

\item[(iii)] As $\e$ tends to $0$, $\Sigma_\e$ converges, away  from  $\{x_1, 
\ldots, x_{k}\}$,  to $\ds \cup_{j=0}^{k} \Pi_j$ in ${\cal C}^\infty$ topology.

\item[(iv)] For all $j=1, \ldots, k$,  the rescaled submanifold  $\e^{-1/n} \, 
(\Sigma_\e -x_j) $  converges to $\alpha_j \, H_{{\cal R}_j}$, where $(\alpha_1,  \ldots,
\alpha_{k}) = \Gamma^{-1} \Lambda$.

\end{enumerate}

\end{theorem}

If $(H1)-(H3)$ are not fulfilled it is still possible to  desingularize 
$\cup_{j=0}^{k} \Pi_j$.  However, this time we only obtain a  minimal submanifold  
which is
not embedded.

\medskip

More generally we should remark that not every configuration  of planes and points can be 
obtained as limit of families of  immersed minimal submanifolds. For example, Ross
\cite{ross}  proved that a minimal two sphere with two punctures minimally  immersed in
${\R}^{4}$ with two simple planar ends has to be  holomorphic, thus giving a severe
restriction on the planes.

\medskip

A number of questions arise naturally form the above result.  In first place it would be 
interesting to know whether conditions  $(H1)-(H3)$ are also necessary for such a
desingularization to exist.  It is easy to check that $(H1)$ is necessary and we believe
that $(H2)$ and $(H2)$ should also be necessary conditions, at least when restricting  the
topological type of the minimal submanifolds. If this turns  out to be the case, since there
are special Lagrangian  configurations of planes which do not satisfy these conditions, it
would mean that they cannot be obtained as  limits of special Lagrangian submanifolds, thus
helping to  understand the possible degenerations of these manifolds.

\medskip

Not unrelated is the problem of determining whether one can  desingularize a 
configuration of special Lagrangian planes through special Lagrangian submanifolds, possibly
leaving free the phase of the calibration to change as in \cite{salur}. In the case of two
planes Lawlor's examples answer obviously the question. A. Brown has recently generalized
this to the case of the  intersection of two special Lagrangian submanifolds with boundary. 
Nothing is known for more than one point of intersection.

\medskip

Finally, we would like to mention the recent work of J. Isenberg, R. Mazzeo and D. 
Pollack  \cite{IMP} where a technical analysis similar to our is developed. 

\section{Comments and examples}

In this section, we give examples of sets of points $x_1, \ldots, x_{k}\in {\R}^n$, 
orthogonal transformations ${\cal R}_1, \ldots, {\cal R}_{k} \in O(n)$, and  matrices ${\cal
A}_0 \in M_n ({\R})$ for which $(H1)$, $(H2)$ and $(H3)$ hold. To simplify the discussion,
let us define 
\[
{\cal S} : = \{(x_1, \ldots, x_{k})\in  {\R}^n \times \ldots \times {\R}^n \, : \,  
x_j \neq x_{j'} \quad \mbox{if}\quad j\neq j'\}.
\]
and 
\[
\Omega : = O(n) \times \ldots \times O(n)
\]

\subsection{Comments on the assumptions}

Observe that 
\[
\mbox{Im} \,  (I - {\cal R}_{j'}^{-1}\, {\cal R}_j) =\mbox{Im} \,  (I - {\cal R}_j^{-1}\, 
{\cal R}_{j'}),
\]
so that condition (H1) is symmetric in $j$ and $j'$. Now, as may easily be checked, 
this assumption guarantees that 
\[
\forall j \neq j', \qquad \qquad \Pi_j \cap \Pi_{j'} =\emptyset ,
\]
while $\Pi_0 \cap \Pi_j =\{x_j\}$ for all $j \geq 0$.

\medskip

Obviously, the set of  $((x_j)_j , ({\cal R}_j )_j) \in  {\cal S} \times  \Omega$ 
such that (H1) is fulfilled is an open set  which is not equal to ${\cal S} \times \Omega$,
since, for example, if  $I - {\cal R}_{j'}^{-1} \, {\cal R}_j$ is invertible for some $j \neq
j'$, then  (H1) does not hold.  Similarly, that the set of  $((x_j)_j , ({\cal R}_j )_j) \in
{\cal S} \times  \Omega$ such that (H2) holds is an open set which is not equal to ${\cal S}
\times \Omega$, since, for example, when all ${\cal R}_j$ are equal, then $\Gamma \equiv 0$
and hence (H2) does not hold.   Finally the set of  $((x_j)_j , ({\cal R}_j )_j, {\cal
A}_0)\in {\cal S} \times \Omega \times  M_n ({\R})$ for which (H3) is fulfilled is also open. 

\medskip

In $M_n ({R})$ we define the relation ${\cal A} \sim {\cal A}'$ if and only if there 
exists ${\cal R}, {\cal R}' \in O(n)$ such that 
\[
{\cal A} = {\cal R} \, {\cal A} \,  {\cal R}'
\]
Now if ${\cal A}_0 \sim {\cal A}'_0$ and if $(H3)$ is fulfilled for ${\cal A}_0$ and
$({\cal R}_j)_j \in \Omega$ then $(H3)$ is fulfilled for   ${\cal A}'_0$ and $({\cal R} \, 
{\cal R}_j\, {\cal R}')_j  \in \Omega$. Hence the set of ${\cal A}\in M_n({\R})$ such that
$(H3)$ is fulfilled for some $((x_j)_j , ({\cal R}_j )_j) \in  {\cal S} \times  \Omega$ only
depends on the coset of ${\cal A}$ in $M_n({\R})/\sim$.

\subsection{Examples}

When $k=2$, we give some examples for which $(H1)$, $(H2)$ and $(H3)$ are fulfilled.
 Without loss of generality, we can assume that $x_1= - x_2 :=e$. Let us assume that both
orthogonal transformations ${\cal R}_i$ leave $e$ unchanged, that is 
\[
{\cal R}_1 \, e = {\cal R}_{2} \, e =e
\]
Therefore 
\[
e \notin \mbox{Im}(I - {\cal R}^{-1}_1 \, {\cal R}_2)
\]
and in particular, $(H1)$ is always fulfilled for such a choice of ${\cal R}_j$. 

\medskip

We would like to compute 
\[
\gamma_{12} =  \frac{1}{2^{n}} \, {\int}_{S^{n-1}} \left( {\cal R}_1  \Theta \cdot 
{\cal R}_2 \Theta  - n \, (\Theta \cdot e)\, (\Theta \cdot e) \right) \, d\theta ,
\]
To this aim, decompose ${\R}^n$ into the direct sum of $E_+ \oplus E_- \oplus 
( \oplus_{i} E_i)$ where $E_\pm$ are the eigenspaces of ${\cal R}_2^{-1}\, {\cal R}_1$
corresponding to the eigenvalues $\pm 1$ and all $E_i$ are two dimensional vector spaces on
which the restriction of ${\cal R}_2^{-1}\, {\cal R}_1$ is a rotation of angle $\theta_i \in
(0, \pi)$.  Since all these spaces are mutually orthogonal, it is easy to compute
\[
\gamma_{12} = -  \frac{\omega_n }{ 2^{n}} \, \left( \frac{2}{n} \, \mbox{dim} 
\, E_{-}  + \frac{2}{n}\,\sum_i (1 -\cos \theta_i ) \right) .
\]
where $\omega_n := |S^{n-1}|$.  Hence the matrix $\Gamma$ is invertible except 
when ${\cal R}_1= {\cal R}_2$. Hence, $(H_2)$ holds whenever ${\cal R}_1 
\neq {\cal R}_2$.
Observe that in this case we have $\gamma_{12} <0$.

\medskip

Finally it remains to check that it is possible to choose ${\cal A}_0 \in M_n ({\R})$ 
such that 
\[
\lambda_j : = -  \int_{S^{n-1}} {\cal A}_0 \Theta \cdot {\cal R}_j \Theta \, 
d \theta  < 0,
\]
for $j=1,2$. To this aim, we define 
\[
L : {\cal A}_0 \in M_n ({\R}) \longrightarrow \Lambda \in {\R}^2 .
\]
This linear map is easily seen to be non zero so, either it is surjective or its image 
is included in a one dimensional space. Assume that the latter is true, 
then there exists
$\alpha, \beta \neq 0$ such that 
\[
\alpha \,  \int_{S^{n-1}} {\cal A}_0 \Theta \cdot {\cal R}_1 \Theta \, d \theta =   
\beta \, \int_{S^{n-1}} {\cal A}_0 \Theta \cdot {\cal R}_2 \Theta \, d \theta ,
\]
for all ${\cal A}_0$. Using this equality for ${\cal A}_0= {\cal R}_1$ and ${\cal A}_0= 
{\cal R}_1$, we conclude that necessarily $\alpha = \pm \beta$ and  
\[
  \int_{S^{n-1}} {\cal R}_1 \Theta \cdot {\cal R}_2 \Theta \, d \theta = \pm \omega_n . 
\]
However, direct computation shows that
\[
\int_{S^{n-1}} {\cal R}_1 \Theta \cdot {\cal R}_2 \Theta \, d \theta = \omega_n   \, 
\left( 1 -  \frac{2}{n} \, \mbox{dim} E_{-} - \frac{2}{n}\,\sum_i (1 -\cos \theta_i ) 
\right) 
\]
Since we already know that this quantity is in absolute value strictly less than 
$|S^{n-1}|$ when ${\cal R}_1 \neq {\cal R}_2$, this implies that $L$ 
is surjective and hence
$L^{-1} (0, +\infty)^2$ is an open nonempty set in $M_n ({\R})$.
 
\medskip

To summarize, we have obtained the 
\begin{lemma}
Assume that $e \in {\R}^n$, ${\cal R}_1, {\cal R}_2 \in O(n)$ are chosen so that 
\[
{\cal R}_1 \, e = {\cal R}_2 \, e =e
\]
then $(H1)$ holds. If in addition ${\cal R}_1 \neq {\cal R}_2$ then $(H_2)$ holds 
and the set of ${\cal A}_0$ such that $(H3)$ holds is an open subset of $M_n({\R})$.
\end{lemma}

\section{Definition of the connection Laplacian on the tangent and normal bundles}

We define all the operators which will be needed in the subsequent sections. We also 
recall some well known properties of these operators.

\subsection{First order differential operators}

To begin with let us define the connections in the tangent bundle and normal bundle of 
a $m$-dimensional submanifold $M$ of ${\R}^n$. Let $V$ be any vector field on a submanifold
$M$ of ${\R}^{N}$ and $e$ a tangent vector field on $M$. We will denote by
$\overline{\nabla}_e V$ the full derivative of $V$ along $e$. 

\medskip

The connection on the tangent bundle $\nabla^{\tau}$, along $e$, applied to the tangent
vector field $T$, is defined to be the orthogonal projection of $\overline{\nabla}_e T$ on
the tangent bundle. We will also write  
\[
\nabla_e^{\tau}T = [ \overline{\nabla}_e T ]^{\tau_M} 
\]
where $[\, \cdot\, ]^{\tau_{M}}$ denotes the orthogonal projection over the tangent 
bundle of $M$.

\medskip

Finally,  we define the connection on the normal bundle $\nabla^{\nu}$, along $e$, 
applied to the normal vector field $V$, is defined to be the orthogonal projection of
$\overline{\nabla}_e N$ on the normal bundle. We will also write  
\[
\nabla_e^{\nu}N = [ \overline{\nabla}_e N ]^{\nu_{M}} 
\]
where $[\, \cdot\, ]^{\nu_{M}}$ denotes the orthogonal projection over the normal 
bundle of $M$.  

\medskip

Let $(e_1, \ldots, e_m)$ be a local orthonormal tangent frame field. For any function 
$f$ defined on $M$, we set
\[
\mbox{grad}_M \, f  : = \sum_{j=1}^{m} \, (e_j f) \, e_j,
\]
 and for any tangent vector field $T$
\[
\mbox{div}_M \, T : = \sum_{j=1}^{m} \, \nabla^{\tau}_{e_j} T \cdot e_j.
\]
where $(e_1, \ldots, , e_m)$ is a local  orthonormal  frame.

\medskip

\subsection{Second order differential operators}

Given $(e_1, \ldots, e_m)$ a local orthonormal tangent frame field, we can define  
$\Delta_M$ the Laplace operator on $M$ acting on the function $f$ by
\[
\Delta_M f  :=  \sum_{j=1}^m \, e_j^2 \, f - \sum_{j=1}^m \, \left( {\nabla^{\tau}_{e_j}e_j} 
\right) \, f.
\]

\medskip

\noindent
The connection Laplacian on the tangent bundle of $M$, acting on the tangent vector 
field $T$, is defined by
\noindent
\[
\Delta^{\tau}_M T  : =  \sum_{j=1}^m  \, \nabla^{\tau}_{e_j}\, \nabla^{\tau}_{e_j} 
T - \sum_{j=1}^m \, \nabla^{\tau}_{\nabla^{\tau}_{e_j}e_j} T
\]
This is just the trace of the invariant second derivative defined by
\[
\left(\nabla^{\tau}\right)^2_{V,W}  : =  \nabla^{\tau}_{V}\, \nabla^{\tau}_{W} - 
\nabla^{\tau}_{\nabla^{\tau}_{V}W}
\]
Let us recall the main properties of $\Delta^\tau_M$. 
\begin{proposition}
\cite{Law-Mic} 
The operator $\Delta^{\tau}_M$ is a negative self-adjoint operator and 
\[
\int_M \Delta^{\tau}_M V \cdot W  = - \int_M \nabla^\tau V \cdot \nabla^\tau W. 
\]
\end{proposition}
Here, by definition
\[
 \int_M \nabla^\tau V \cdot \nabla^\tau W := \sum_{j=1}^m \int_M \nabla_{e_j}^\tau V 
\cdot \nabla^\tau_{e_j} W. 
\]

\medskip

\noindent
Finally, the connection Laplacian operator on the normal bundle of $M$, acting on the 
normal vector field $N$, is defined to be 
\[
\Delta^{\nu}_M N : = \sum_{j=1}^m \, \nabla^{\nu}_{e_j}\, \nabla^{\nu}_{e_j} - \sum_{j=1}^{m} 
\, \nabla^{\nu}_{\nabla^{\tau}_{e_j}e_j} \, N
\]
This is just the trace of the invariant second derivative defined by
\[
\left(\nabla^{\nu}\right)^2_{V,W}  : =  \nabla^{\nu}_{V}\, \nabla^{\nu}_{W} - \nabla^{\nu}_{
\nabla^{\tau}_{V}W}
\]
Let us recall the main properties of $\Delta^\nu_M$. 
\begin{proposition}
\cite{Law-B} 
The operator $\Delta^{\nu}_M$ is a negative self-adjoint operator and 
\[
\int_M \Delta^{\nu}_M V \cdot W  = - \int_M \nabla^\nu V \cdot \nabla^\nu W. 
\]
\end{proposition}
Here, by definition
\[
\int_M \nabla^\nu V \cdot \nabla^\nu W := \sum_{j=1}^n \int_M \nabla_{e_j}^\nu V \cdot 
\nabla^\nu_{e_j} W. 
\]
\medskip

\noindent
Observe that all these definitions do not depend on the choice of the local orthonormal 
tangent frame field $(e_1, \ldots, e_m)$. 

\subsection{Operators of order zero}

The second fundamental form is a section of the bundle $T^\star (M) \otimes T^\star (M) 
\otimes N (M)$ defined by 
\[
B_{V,W} : = \nabla^{\nu}_V W
\]
In other words, at any point $p \in M$, $B_p$ represents a symmetric bilinear map from 
$T_p (M)$ into $N_p(M)$. 

\medskip

Given any local orthonormal tangent frame field $(e_1, \ldots, e_m)$ we define
\[
B_{i,j} = B_{e_i,e_j} =\nabla^\nu_{e_i}e_j
\]
This allows to define the linear operator acting on normal vector fields
\[
{\cal B} \,  (N) : = \sum_{i,j=1}^m (B_{ij} \cdot N ) \, B_{i,j} = \sum_{i,j=1}^m 
(\nabla^\nu_{e_i}e_j \cdot N ) \, \nabla_{e_i}e_j
\]

\section{Differential forms}

It will be useful to translate some of the previously defined operator in the 
language of differential forms. Let $\Omega^p (M)$ denote the space of $p$-forms on $M$. We
denote by $d^p$ the exterior derivative
\[
d^p : \Omega^p (M) \longrightarrow \Omega^{p+1}(M)
\]
and define the operator 
\[
\delta^p : \Omega^{p} (M) \longrightarrow \Omega^{p-1}(M)
\]
 by
\[
\delta^p = (-1)^{d(p+1)+1}\, \star \, d \, \star
\]
where $\star$ denotes the usual Hodge operator. 

\medskip

The Hodge Laplacian on $p$ forms is defined by
\[
\Delta^p = - \left( d^{p-1}\, \delta^p + \delta^{p+1} \, d^p \right).
\]
Observe that, with this definition, $\Delta^p$ is a negative operator.

\medskip

The inner product on $\Omega^p(M)$ is defined by
\[
<\omega, \tilde{\omega}>_{\Omega^p(M)}  : = \int_M \omega \wedge \star \, \tilde{\omega}
\]
Recall that, granted the definition of $\delta^p$, for any $p-1$-form $\omega$ and any 
$p$-form $\tilde{\omega}$, we have 
\[
\int_M d^p \omega \wedge \star \, \tilde{\omega} = \int_M  \omega \wedge \star \, \delta^p 
\, \tilde{\omega}  
\]
and hence, for any $p$-forms $\omega, \tilde{\omega}$, we have
\[
\int_M \Delta^p \omega \wedge \star \, \tilde{\omega} = - \int_M  \left( d^{p} \omega 
\wedge \star \, d^p \tilde{\omega}  +\delta^{p} \omega \wedge \star \delta^p \tilde{\omega} 
\right)
\]

\medskip

We will use the Hodge decomposition Theorem we recall now
\begin{theorem}
Let $M$ be a compact submanifold of ${\R}^n$. Then $\Omega^p (M)$ can be uniquely 
decomposed as
\[
\Omega^p (M) = d^{p-1}\, \Omega^{p-1}(M)\oplus \delta^{p+1}\, \Omega^{p+1}(M)\oplus 
\mbox{Har}\, ^p(M)
\]
where $\mbox{Har}\, ^p(M)$ is the set of harmonic $p$-forms.
\end{theorem}

\medskip

One can identify functions on $M$ with $0$-forms in the obvious way and, using the metric 
on $M$, one can also identify tangent vector fields with $1$-forms in the following way.
Assume that we have chosen a  local orthonormal tangent frame field $(e_1, \ldots, e_m)$ on
$M$ and that the metric is then given by $g= (g_{ij})_{i,j}$ so that 
\[
\e_j =  \frac{1}{\sqrt{g_{jj}}} \, \frac{\del}{\del x_j}
\]
Then to any vector field
\[
T = \sum_{j=1}^m T_i \, e_i
\]
we can associate the $1$-form
\[
\omega = \sum_{j=1}^m \sqrt{g_{jj}} \, T_j \, dx_j
\]
This identification is coherent with the definition of the inner product on $\Omega^1(M)$ 
and on $T(M)$, namely, if the vector field $T$ is associated to the $1$-form $\omega$, we
have 
\[
\|T\|^2_{L^2(M)}:= \int_M |T|^2 = \int_{M} \omega \,\wedge \,  \star \, \omega : = 
\|\omega\|^2
\]

\medskip

Granted this identification of vector fields with $1$-forms, we can identify $\Delta^1$, 
the Hodge Laplacian on $1$-forms, with some symmetric second order differential operator
acting on tangent vector fields. We will still denote by $\Delta^1$ this operator which is
now defined on tangent vector fields. The relation between this operator and the connection
Laplacian on the tangent bundle which we have defined in the previous section is given by the
following 
\begin{theorem}
\cite{Law-Mic}
The difference between the Hodge Laplacian on $1$-forms and the connection Laplace operator 
on the tangent bundle is given by
\[
\Delta^1 = \Delta^\tau_M - \mbox{Ric},
\]
where $\mbox{Ric}$ denotes the Ricci tensor.
\end{theorem}

\section{The hyperbola}

The hyperbola $H_I$ in dimension $2 \, n$ is parameterized by
\[
X (s, \theta) : = \left( \frac{\cos s}{(\sin (ns))^{1/n}} \, \Theta (\theta) ,  
\frac{\sin s}{(\sin (ns))^{1/n}} \, \Theta (\theta) \right)
\]
where $s \in (0, \frac{\pi}{n})$ and where 
\[
\theta \longrightarrow \Theta (\theta) \in S^{n-1}
\]
is a parameterization of the $n-1$ dimensional sphere. It will be convenient to assume 
that this parameterization is chosen in such a way that 
\begin{equation}
\del_{\theta_j} \Theta \cdot  \del_{\theta_k} \Theta =0 \qquad \mbox{if} \qquad j \neq k .
\label{eq:norm-param}
\end{equation}

\medskip

\begin{remark}
When $n=2$, we have the equivalent definition of the hyperbola given by
\[
z \in {\C}\setminus \{0\}\longrightarrow  \left( z, \frac{1}{\bar{z}} \right)  \in {\C}^2
\]
\end{remark}

\medskip

For notational convenient , we will frequently write $S$ instead of $S^{n-1}$ in 
subscripts or superscripts. It will also be convenient to identify ${\R}^{2n}$ with 
${\C}^n$
using
\[
{\R}^n \times {\R}^n \ni  (x, y)  \sim x+iy \in {\C}^n .
\]
With this identification we will write 
\[
X (s, \theta) : = \frac{e^{is}}{(\sin (ns))^{1/n}} \, \Theta (\theta).
\]

\medskip

{\bf The tangent bundle}. The tangent space of $H_I$ is spanned by the following set 
of vectors
\[
\del_s X = - \frac{e^{i(1-n)s}}{(\sin (ns))^{1+1/n}}\, \Theta
\]
and, for all $j = 1, \ldots, n-1$
\[
\del_{\theta_j}  X  =  \frac{e^{i s}}{(\sin (ns))^{1/n}}\, \del_{\theta_j} \, \Theta
\]

\noindent
In order to simplify the notations, we will write for short 
\[
\Theta_j := \del_{\theta_j} \, \Theta.
\]
and we will define the vectors
\[
e_0 := e^{i(1-n)s}\, \Theta \qquad \mbox{and}\qquad  e_j : = e^{is}\, \Theta_j
\]
so that, thanks to (\ref{eq:norm-param}), $(e_0, \ldots, e_{n-1})$ is an orthonormal basis of  $T(H_I)$ at the point $X$.

\medskip

\noindent
Finally, we define, for all $j=1, \ldots, n-1$
\[
\e_j := \frac{\Theta_j}{|\Theta_j|}
\]
so that, $(\e_1, \ldots, \e_{n-1})$ is an orthonormal basis of $T (S^{n-1})$ at the 
point $\Theta$.

\medskip

{\bf The normal bundle}. The normal space to $H_I$ is spanned by the following vectors 
\[
N_0 : = i \, e^{i(1-n)s}\, \Theta,
\]
and, for $j=1, \ldots, n-1$
\[
N_j : = i  \, e^{is}  \, \frac{\Theta_j}{|\Theta_j|}.
\]
Hence any normal vector field on $H_I$ can be written as
\[
V :=  i \, e^{i(1-n)s}\, f \, \Theta + i \, e^{is}\, T,
\]
where $f$ is a real valued function and, for all $s \in (0, \frac{\pi}{n})$, $T(s, \cdot)$ 
is a tangent vector field on $S^{n-1}$. 

\section{Expansion of the lower end of the hyperbola}

Recall that the hyperbola is parameterized by
\[
X(s, \theta) = \frac{e^{is}}{(\sin (ns))^{1/n}}\, \Theta
\]
We set 
\[
r := \frac{\cos s}{(\sin (ns))^{1/n}}
\]
We can expand $r$ in terms of $s$ as $s$ tends to $0$. We get explicitely
\[
r = (ns)^{-1/n}\, \left( 1 + {\cal O} (s^2)\right)
\]
Which in turn yields an expansion of $s$ in terms of $r$, as $s$ tends to $0$
\[
s = \frac{r^{-n}}{n}\, \left( 1  + {\cal O} (r^{-2n})\right)
\]
Finally, we get the expansion as $s$ tends to $0$
\[
\frac{\sin s}{(\sin (ns))^{1/n}} = \frac{r^{-n}}{n} \, \left( 1 + {\cal O}(r^{-2n})\right)
\]
Hence the hyperbola is parameterized by
\[
X (r, \theta) =  r \Theta + i \,  \frac{r^{1-n}}{n} \, \left( 1 + {\cal O}(r^{-2n})\right) \, \Theta
\]
as $r$ tends to $+\infty$.

\medskip

We now consider the Hyperbola scaled by a factor $(n \, \beta \, \e )^{1/n}$ for some 
$\beta >0$ and $\e >0$, namely
\[
X_\beta (s, \theta) = (n \, \beta \, \e )^{1/n}  \, \frac{e^{is}}{(\sin (ns))^{1/n}}\, 
\Theta
\]
We set $\rho := (\beta \, \e)^{\frac{1}{n}} \, r$ and, as $s$ tends to $0$ we can also 
write
\begin{equation}
X_\beta (\rho, \theta) = \rho \Theta + i \,  \e \, \beta \, \rho^{1-n} \, \left( 1 + 
{\cal O}(\e^{2} \, \rho^{-2n})\right) \, \Theta .
\label{eq:expansion}
\end{equation}

\section{The Linearized mean curvature operator about the hyperbola}

In order to compute the linearized mean curvature operator, we collect the results of 
the Appendix 1.  We recall from \cite{Law-B} that the linearized mean curvature operator is
given by
\[
L_{H} = \Delta^{\nu}_H + {\cal B}
\]
In the local chart given by $X$, we have explicitely
\begin{proposition}
The linearized mean curvature operator about the hyperbola $H_I$ reads
\[
\begin{array}{rcllll}
(\sin (ns))^{-\frac{2}{n}} \, L_H \,  V  & =   &   i \,  e^{i(1-n)s} \, \left[  
(\sin (ns))^{2-\frac{2}{n}} \, \del_s \left( (\sin (ns))^{\frac{2}{n}} \, \del_s f \right) +
\Delta_S f   - (n-1)  f   \right. \\[3mm]
                                                               &      &	
\qquad \qquad                                    \left. + (n^2-1) \, \sin^2 (ns)  
\, f  - 2
\, \cos (ns) \, \mbox{\em div}_S \,  T \right]  \, \Theta\\[3mm]
                                                               &  +  &   i 
\, e^{is} \,  \left[ (\sin (ns))^{2-\frac{2}{n}} \, \del_s \left( (\sin (ns))^{2/n} \, \del_s
T \right) +  \Delta^{\tau}_S T  -  T \right. \\[3mm]
                                                               &      &	
\qquad \qquad 	        	           \left.  + 3 \,  \sin^2 (ns)\, T \, +  2 \,  \cos (ns) 
\, 
\mbox{ \em grad}_S \,  f  \right]
\end{array}
\]
\end{proposition}

\medskip

In order to study the spectral properties of $L_H$, it will be useful to identify 
tangent vector fields on  $S^{n-1}$ with $1$-forms on $S^{n-1}$. This identification yields a
natural identification of normal vector fields on $H_I$
\[
V : = f \, i \, e^{i(1-n)s}\, \Theta + i \, e^{is}\, T  \in N(H)
\]
with 
\[
U : = (f , v ) \in C^{\infty} ({\R}\times S^{n-1} ; \Omega^0  \times \Omega^1 )
\]
where $v \, (s, \cdot)$ is the $1$-form corresponding to $T(s, \cdot)$.

\medskip

Granted this identification we can identify $L_H$ with the following linear operator
\[
\begin{array}{rlll}
(\sin (ns))^{-\frac{2}{n}} \, L_H \,  (f, v)    &  =  &   (\sin (ns))^{2-\frac{2}{n}} 
\, \del_s \left( (\sin (ns))^{\frac{2}{n}} \, (\del_s f, \del_s v) \right) + (\Delta^0 f ,
\Delta^1 v) \\[3mm]
                                                                      &  + 
 & 2 \,  \cos (ns) \, (\delta^1 \,  v , d^0 \, f)   +   ((1-n) \, f , (n-3) \, v)   \\[3mm]
                                                                      &  + 
 &  \sin^2 (ns) \, ((n^2-1) f , 3 v)  
\end{array}
\]

\medskip

Also, in order to have a better understanding of the structure of this operator we 
define the variable $t$ by
\[
dt  : = \frac{1}{\sin (ns)} \, ds,
\]
with $t (\frac{\pi}{2n}) = 0$. We obtain explicitely that 
\[
e^{-nt} = \frac{\sin (ns)}{1- \cos (ns)}.
\]
which implies that 
\[
\sin (ns) = ( \cosh (nt) )^{-1} \qquad \mbox{and}\qquad \cos (ns) = - \tanh (nt).
\]
Then one checks directly that the operator
\[
D :=  (\sin (ns))^{2} \, \del_s \left( (\sin (ns))^{\frac{2}{n}} \, \del_s \right) ,
\]
becomes
\[
(\sin (ns))^{- \frac{n+2}{2n}}  \, D  \, (\sin (ns))^{\frac{n-2}{2n}}  = \del_t^2  - 
\frac{(n-2)^2}{4} - \frac{n^2-4}{4} \, \frac{1}{\cosh^2 (nt)}  .
\]

\medskip

Hence the study of $L_H$ is equivalent to the study of the conjugate operator
\[
\begin{array}{rlll}
{\cal L}_H \,  (f, v)   & = & \ds   (\del_t^2 f, \del_t^2 v)  + (\Delta^0 f , \Delta^1 v)  
- 2 \,  \tanh (nt) \, (\delta^1 \,  v , d^0 \, f)   \\[3mm]
                                              & - &  \ds   \left( \frac{n^2}{4} \, f , 
\frac{(n-4)^2}{4} \, v \right) + \frac{1}{\cosh^2 (nt)}  \, \left( \frac{3 n^2}{4} \, f,
\frac{16-n^2}{4} \,  v\right)   
\end{array}
\]

Observe that, when $t$ tends to $- \infty$ the operator ${\cal L}_H$ is equivalent to the 
following differential operator
\[
{\cal L}_0 \,  (f, v)   = \ds   (\del_t^2 f, \del_t^2 v)  + (\Delta^0 f , \Delta^1 v)  + 2 \,   (\delta^1 \,  v , d^0 \, f)    -  \left( \frac{n^2}{4} \, f , \frac{(n-4)^2}{4} \, v \right) 
\]

\section{Eigendata of $\Delta^0$ and $\Delta^1$ on $S^{n-1}$}

The spectrum of $\Delta^0$ on $S^{n-1}$ is well known and given by
\[
\sigma(\Delta^0) =\{ k \, (n-2+k) \, | \, k \geq 0\}
\]

\medskip

The spectrum of $\Delta^1$ is also well known. In dimension $n=2$ or $n=3$ this spectrum is 
simply given by 
\[
\sigma (\Delta^1) =\{ k \, (n-2+k) \, | \, k \geq 1\}
\]
and all $1$-eigenforms of $\Delta^1$ are image of eigenfunctions of $\Delta^0$ by the operator $d^0$.

\medskip

In dimension $n \geq 4$, things are slightly more involved. The spectrum of $\Delta^1$ 
can be decomposed into two disjoint subsets
\[
\sigma (\Delta^1) = \sigma_{ex} (\Delta^1) \cup \sigma_{coex}(\Delta^1)
\]
where 
\[
\sigma_{ex} (\Delta^1) =\{ k \, (n-2+k) \, | \, k \geq 1\}
\]
corresponds of the eigenvalues associated to exact $1$-eigenforms (namely $1$-eigenforms 
which belong to the image of $\Omega^0(S^{n-1})$ by  $d^0$) and where
\[
\sigma_{coex} (\Delta^1) =\{ (k+1) \, (n-3+k) \, | \, k \geq 1\}
\]
corresponds of the eigenvalues associated to coexact $1$-eigenforms  
(namely $1$-eigenforms which belong to the image of $\Omega^2(S^{n-1})$ by $\delta^2$).  It
may be observed that all exact $1$-eigenforms are image of eigenfunctions of $\Delta^0$ by
$d^0$ and all coexact $1$-eigenforms are image of exact $2$-eigenforms of $\Delta^2$ by
$\delta^2$.

\medskip

We define 
\[
{\cal V}^0 := \left\{ (f, 0) \, : \, d^0\,  f=0  \right\}
\]
Next, for all $k\geq 1$
\[
{\cal V}^k_{ex} := \left\{ (f, v ) \, : \, (\Delta^0 f , \Delta^1 v) = - k(n-2+k) \, (f, v) 
\quad \mbox{and}\quad d^1 v =0 \right\}
\]
the eigenspace corresponding to the eigenvalue $k(n-2+k)$ and to exact $1$-forms, 
and 
\[
{\cal V}^k_{coex} := \left\{ (0, v) \, : \, \Delta^1  v = - (k+1)(n-3+k) \, v \quad 
\mbox{and}\quad \delta^1 v=0 \right\}
\]
the eigenspace corresponding to the eigenvalue $k(n-2+k)$ and to coexact $1$-forms.

\medskip

This decomposition corresponds to the decomposition for $1$-forms on $S^{n-1}$ given 
by Hodge's decomposition Theorem
\[
\Omega^1 (S^{n-1}) = d^0 \, \Omega^0 (S^{n-1}) \oplus \delta^2 \, \Omega^2 (S^{n-1}),
\]
since there is no harmonic $1$-form on $S^{n-1}$. The set of exact $1$-eigenforms is total 
in $d^0 \, \Omega^0 (S^{n-1})$ and the set of coexact $1$-eigenforms is total in $\delta^1 \,
\Omega^2(S^{n-1})$. In addition the two spaces $d^0 \, \Omega^0 (S^{n-1})$ and $\delta^2 \,
\Omega^2 (S^{n-1})$ are orthogonal.

\medskip

\begin{remark}
This decomposition of the space of $1$-forms induces a natural decomposition of the space 
of tangent vector fields into vector fields corresponding to $d^0 \, \Omega^0 (S^{n-1}) $ and
vector fields corresponding to $ \delta^2 \, \Omega^2 (S^{n-1})$. Moreover these two sets of
vector fields are $L^2$ orthogonal.
\end{remark}

\section{Indicial roots}

To begin with let us compute the indicial roots corresponding to ${\cal L}_H$. First, 
observe that, in dimension $n\geq 4$, we can decompose any $1$-form $v= v_{ex}+ v_{coex}$
into the sum of an exact $1$-form on $S^{n-1}$ and a coexact $1$-form on $S^{n-1}$. In
dimension $n=2$ and $n=3$, we set $v_{coex}=0$. Now, if $(f, v)$ is a solution of the
homogeneous equation ${\cal L}_H (f, v)=0$ then $v_{coex}$ satisfies
\[
\begin{array}{rlll}
{\cal L}_H \,  (0 , v_{coex})  &  =  &   \ds  (0, \del_t^2 v_{coex})  + (0, 
\Delta^1 v_{coex})  - \left( 0, \frac{(n-4)^2}{4} \, v_{coex} \right)  \\[3mm]
				&  +  &  \ds  \frac{1}{\cosh^2 (nt)}  \, \left( 0, \frac{16-n^2}{4} \,  
v_{coex} \right)    
\end{array}
\]
Recall that, restricted to coexact forms, the spectrum of $\Delta^1$ is given 
by the set of $(k+1) \,  (n-3+k)$, for $k\geq 1$. In order to find the indicial roots
corresponding to this operator, we look for solutions of the homogeneous problem of the form
\begin{equation}
v_{coex} = a \, \psi
\end{equation}
where $\psi$ is a coexact eigenform of $\Delta^1$ corresponding to the eigenvalue 
$(k+1) \, (n-3+k)$ and where $a$ is a scalar function only depending on $t$. It is easy to
see that the scalar function $a$ is a solution of the following ordinary differential 
equation
\begin{equation}
\ddot{a} -\left( \frac{n-2}{2} + k \right)^2 \, a +  \frac{16-n^2}{4} \, 
\frac{a}{\cosh^2 (nt)} =0
\end{equation}
The behavior of the solutions when $t$ tends to either $\pm \infty$ is given by  
$e^{\gamma_k^{\pm} \, t}$ where 
\begin{equation}
\gamma_k^\pm = \pm \left( \frac{n-2}{2} + k \right) 
\end{equation}
These are the indicial roots of ${\cal L}_H$ when this operator is restricted to 
coexact $1$-forms.

\medskip

We now turn to the study of indicial roots corresponding to $\omega_{ex}$. 
Observe that, this time
\begin{equation}
\begin{array}{rlll}
{\cal L}_H \,  (f, v_{ex})   & = &    \ds (\del_t^2 f, \del_t^2 v_{ex})  + 
(\Delta^0 f , \Delta^1 v_{ex})   - 2 \,  \tanh (nt) \, (\delta^1 \,  v_{ex} , d^0 \, f) 
\\[3mm] &  - &  \ds   \left( \frac{n^2}{4} \, f , \frac{(n-4)^2}{4} \, v_{ex} \right)  + 
\frac{1}{\cosh^2 (nt)}  \, \left( \frac{3 n^2}{4} \, f, \frac{16-n^2}{4} \, 
v_{ex}\right)      
\end{array}
\label{eq:ind-root}
\end{equation}
Assume that $\phi$ is an eigenvalue of $\Delta^0$ associated to the eigenvalue 
$k \, (n-2+k)$. Then, we look for solutions of (\ref{eq:ind-root}) of the form
\[
f= a \, \phi \qquad \mbox{and}\qquad v_{ex} = b \, d^0 \phi
\]
where $a$ and $b$ are scalar functions only depending on $t$. We obtain the following 
system of ordinary differential equation
\[
\ddot{a} - \frac{n^2}{4}\, a  + \frac{3 n^2}{4} \,  \frac{a}{\cosh^2 (nt)} =  0  
\]
when $k=0$ and 
\[
\left\{
\begin{array}{lllll}
\ds \ddot{a} -k \,  (n-2+k) \, a - \frac{n^2}{4}\, a - 2 \, \tanh (nt) \, k \, 
(n-2+k) \, b  +   \frac{3 n^2}{4} \,  \frac{a}{\cosh^2 (nt)} & = & 0 \\[3mm] 
\ds \ddot{b} -k \, (n-2+k) \, b - \frac{(n-4)^2}{4}\, b - 2 \, \tanh (nt) \, a  + 
 \frac{16 - n^2}{4} \,  \frac{b}{\cosh^2 (nt)} & = & 0 \\[3mm] 
\end{array}
\right.
\]
for $k\neq 0$. With little work one finds for $k\neq 0$ the asymptotic behavior of 
$a$ and $b$  at both $\pm \infty$ is governed by the following sets of indicial roots
\[
\mu^{\pm}_k=\pm \left( \frac{n}{2} +k\right) \qquad \mbox{and}\qquad \nu^{\pm}_k=\pm 
\left( \frac{n-4}{2} +k\right) 
\]
and for $k=0$, we find 
\[
\mu^{\pm}_0= \pm \frac{n}{2}.
\]

It is worth mentioning that the operators ${\cal L}_H$ and ${\cal L}_0$ have the same 
indicial roots.

\section{The Jacobi fields}

Some Jacobi fields of $L_H$ are very easy to obtain since they correspond to  geometric 
transformations of the hyperbola. In this paragraph we consider $L_H$ instead of ${\cal L}_H$
since these Jacobi fields are easier to describe for $L_H$.

\medskip

\noindent
{\bf Jacobi fields corresponding to translations} Let $a, b \in {\R}^n$ be given, the 
Jacobi field corresponding to the translation of vector $e^{i\alpha}\, a$, namely
\[
x+ i \, y \in {\C}^n\longrightarrow (x+\cos \alpha \, a)+ i \, (y+\sin \alpha \, b)\in 
{\C}^n
\]
is the projection on the normal bundle of the constant vector $e^{i\alpha}\, a$. We 
obtain explicitely
\begin{equation}
\Phi_t (\alpha, a)  = i \, e^{i(1-n)s}\, \sin ((n-1)s + \alpha) \, (a\cdot \Theta ) \, 
\Theta + i \, e^{is} \,  \sin ( \alpha -s) \, (a - (a\cdot \Theta) \, \Theta)  
\end{equation}

\medskip

\noindent
{\bf Jacobi field corresponding to a dilation} This Jacobi field corresponding to 
the dilation
\[
x+i \, y \in {\C}^n \longrightarrow (1+ \delta ) \, (x+i \, y)\in {\C}^n
\]
is obtained by projecting over the normal bundle the infinitesimal dilation
\begin{equation}
\delta \, \frac{e^{i s}}{(\sin (ns))^{\frac{1}{n}}}\, \Theta
\end{equation}
We obtain 
\begin{equation}
\Phi_d (\delta) = \delta \,  ( \sin (ns))^{1 - \frac{1}{n}} \, i \, e^{i(1-n)s}\, \Theta
\end{equation}

\medskip

\noindent
{\bf Jacobi fields corresponding to the action of $SU(n)$} Let $A$ be some 
$n \times n$ symmetric matrix. The Jacobi fields corresponding to the action of
\[
x+i \, y \in {\C}^n \longrightarrow e^{iA} \, (x+i \, y)\in {\C}^n
\]
are obtained by projecting over the normal bundle the infinitesimal action
\[
i  \, \frac{e^{is}}{(\sin (ns))^{1/n}}\, A \, \Theta
\]
We obtain explicitely
\begin{equation}
\Phi_{SU(n)}(A) = (\sin (ns))^{-\frac{1}{n}}  \left( i \, e^{i (1-n)s}  \, 
\cos (ns)\, (A \, \Theta \cdot \Theta)\, \Theta + i \, e^{is}
\,  (A \, \Theta - (A\, \Theta \cdot \Theta) \, \Theta) \right).
\end{equation}

\medskip

\noindent
{\bf Jacobi fields corresponding to the action of $O(2n)/ SU(n)$} Let $A$ be some 
$n \times n$ antisymmetric matrix. The Jacobi fields corresponding to the action of
\[
x+i \, y \in {\C}^n \longrightarrow (e^{-A} \, x+i \, e^{A} \, y)\in {\C}^n
\]
are obtained by projecting over the normal bundle the infinitesimal action
\[
 \frac{e^{- is}}{(\sin (ns))^{\frac{1}{n}}}\, A \, \Theta
\]
We obtain explicitely
\begin{equation}
\Phi_{O(2n)}(A) =  (\sin (ns))^{-\frac{1}{n}} \, \sin (2s) \,  i \, e^{is}\,  A 
\, \Theta .
\end{equation}

The Jacobi fields corresponding to the action of
\[
x+i \, y \in {\C}^n \longrightarrow (\cosh A \, (x+i \, y) + \sinh A \, (y+i \, x))\in {\C}^n
\]
are obtained by projecting over the normal bundle the infinitesimal action
\[
i \,  \frac{e^{-is}}{(\sin (ns))^{\frac{1}{n}}}\, A \, \Theta
\]
We obtain explicitely
\[
\Psi_{O(2n)}(A) =  (\sin (ns))^{-\frac{1}{n}} \, \cos (2s) \,  i \, e^{is}\,  A \, \Theta .
\]

\section{The maximum principle}

We want to prove the following result
\begin{proposition}
Assume that $n \geq 4$ and assume that $U$ is a solution of ${\cal L}_H \, U =0$ in 
$(t_1, t_2) \times S^{n-1}$, (with $U =0$ at the boundary if $t_1>- \infty$ or $t_2 < +
\infty$). Assume that 
\[
|U|\leq (\cosh t)^{-\nu},
\]
for some $\nu > \frac{2-n}{2}$. Further assume that, for every $t \in (t_1, t_2)$,  
$U(t, \cdot)$ is orthogonal to ${\cal V}^0$  in the $L^2$-sense on $S^{n-1}$. Then $U \equiv
0$.
\label{pr:4}
\end{proposition}

\noindent
{\bf Proof :} We decompose $U =(f, v)$ into $v = v_{ex}+ v_{coex}$ where $v_{ex}$ 
is an exact $1$-form on $S^{n-1}$ and where $v_{coex}$ is a coexact $1$-form on $S^{n-1}$. 

\medskip

{\bf Step 1.}  We multiply the equation ${\cal L}_H (f, v)$ by $(0, v_{coex})$ and 
integrate over $(t_1, t_2) \times S^{n-1}$. We obtain
\[
0  = \ds  \int  |\del_t v_{coex}|^2 + \int |d^1 v_{coex}|^2  +  \ds \frac{n^2}{4} 
\, \int  |v_{coex}|^2 -  \frac{n^2-16}{4}\, \int (\cosh (nt))^{-2}\, |v_{coex}|^2
\]
which already implies that $v_{coex}\equiv 0$ since $n \geq 4$. It therefore remains 
to prove that $(f, v_{ex})\equiv 0$.

\medskip

The proof is now quite involved and, in order to simplify the notations, we set
\[
\begin{array}{lllllllllll}
A :=  \ds \int  f^2,      & B  := \ds  \int (\cosh (nt))^{-2}\, f^2, &  C  := \ds \int  
|\del_t f|^2 , &  D  := \ds \int  |d^0 \,  f|^2  \\[3mm]
A' :=  \ds \int  |v_{ex}|^2 ,  &  B'  :=  \ds \int (\cosh (nt))^{-2}\, |v_{ex}|^2, &  C'  
:=  \ds \int |\del_t v_{ex}|^2, &  D'  :=  \ds \int |\delta^1 v_{ex}|^2
\end{array}
\]
and 
\[
E : = \int  \tanh (nt) \, d^0\, f \wedge \star \, v_{ex} = -  \int \tanh (nt) \,  f \, 
\delta^1 \,  v_{ex})
\]

\medskip

{\bf Step 2.} We multiply the equation ${\cal L}_H (f, v)$ by $(f, v_{ex})$ and integrate. 
This time, we obtain
\[
\begin{array}{llll}
 \ds  \int ( | \del_t f |^2 +  |\del_t v_{ex}|^2) + \int ( |d^0 f|^2 + |\delta^1 v_{ex}|^2)
+  \ds \int \left(  \frac{n^2}{4} \,  f^2 + \frac{(n-4)^2}{4} \, |v_{ex}|^2\right) \\[3mm]
 =  \ds \int (\cosh (nt))^{-2}\, \left(  \frac{3n^2}{4}\,  f^2+ \frac{16-n^2}{4}\, |v_{ex}|^2
\right) -  2 \, \ds \int  \tanh (nt) \, (d^0\, f \wedge \star \, v_{ex} + f \, \delta^1 \, 
v_{ex})
\end{array}
\]
With our notations this can be written as
\begin{equation}
\ds  C+ C' + D+ D' +  \frac{n^2}{4} A + \frac{(n-4)^2}{4} A'  - \frac{3 n^2}{4} B + 
\frac{n^2-16}{4}B'  = -  4 \, E
\label{eq:est-1}
\end{equation}

\medskip

Using with Cauchy-Schwarz inequality and integration by parts, we estimate
\[
\ds  \left| \ds \int  \tanh (nt)  \, f \, \delta^1 \, v_{ex} \right|   \leq  \ds 
\left( \int |\delta^1 v_{ex}|^2 \right)^{\frac{1}{2}}  \left(\int f^2 - \int (\cosh
(nt))^{-2} \, f^2\right)^{\frac{1}{2}}.
\]
This inequality can also be written as 
\begin{equation}
E^2 \leq D' \, (A-B).
\label{eq:est-2}
\end{equation}

\medskip

Now, we use the fact that, for all $t$, $f(t, \cdot)$ is $L^2$-orthogonal to the 
first eigenfunction on $S^{n-1}$ hence
\[
\int_{S^{n-1}} |d^0 \, f |^2 \geq (n-1) \, \int_{S^{n-1}} f^2
\]
Integrating over $t$ yields
\[
\int  |d^0 \, f|^2 \geq (n-1) \, \int \, f^2
\]
Otherwize stated 
\begin{equation}
D \geq (n-1) \, A.
\label{eq:est-3}
\end{equation}

\medskip

Finally, we use an integration by parts to prove that
\[
\ds n \, \int (\cosh (nt))^{-2}\, f^2  =   - \ds 2 \, \int \tanh (nt) \, f \, \del_t f
\]
Using Cauchy-Schwarz inequality, we conclude that
\[
\ds n \, \int (\cosh (nt))^{-2}\, f^2 \leq  \ds 2 \, \left( \int |\del_t f|^2\right)^{
\frac{1}{2}}  \left( \int f^2 - \int (\cosh (nt))^{-2}\, f^2\right)^{\frac{1}{2}}
\]
A similar inequality holds with $f$ replaced by $v_{ex}$. This yields 
\begin{equation}
n^2 \, B^2 \leq 4 \,  C \, (A-B)^2 \qquad \qquad n^2 \, {B'}^2 \leq 4 \,  C' \, (A'-B')^2 
\label{eq:est-4}
\end{equation}

\medskip

{\bf Step 3.} Assume that $A\neq B$ and $A'\neq B'$, we can collect the previous inequalities  (\ref{eq:est-2})-(\ref{eq:est-4}) to eliminate $C$, $C'$, $D$ and $D'$ in (\ref{eq:est-1}). With little work we find 
\[
\ds  \frac{n^2}{4} \frac{(2B -A)^2}{A-B}+ \ds  \frac{1}{4} \frac{(4B'+(n-4)A') ^2}
{A'-B'} + (n-5) \, A +  4\, B+  \frac{(E -A+B)^2}{A-B}\leq 0
\]
It is now an easy exercice to check that, when $n\geq 4$, the above inequality 
implies that $A=A'=0$. hence $f=0$ and $v_{ex}=0$.

\medskip

If  $A=B$, we first use (\ref{eq:est-4} to conclude that $B=0$,  hence  we conclude 
that $f=0$. In this case, it readily follows from (\ref{eq:est-1}) that $A'=0$ which implies
that $v_{ex}=0$.

\medskip

If  $A'=B'$, we first use (\ref{eq:est-4} to conclude that $B'=0$, so $v_{ex}=0$. 
In this case,  (\ref{eq:est-1}) reduces to 
\[
C+D + \frac{n^2}{4} \, A - \frac{3n^2}{4} \, B =0
\]
And, using (\ref{eq:est-2})-(\ref{eq:est-4}) we conclude that $A=0$ which implies 
that $f =0$.

\medskip

The proof is therefore complete. \hfill $\Box$

\medskip

When $n=3$ some modifications are needed in the corresponding statement. Indeed, we 
shall prove, and this will be sufficient for our purposes that the maximum principle as
stated in Proposition~\ref{pr:4} holds on any interval containing a large interval centered
at the origin. This additional assumption is probably not needed but we have not been able to
get ride of it. Let us emphasize that the result needed below is just what we need for all
the remaining analysis to hold.
\begin{proposition}
Assume that $n = 3$. Then there exists $\hat t_0 >0$ such that, if $U$ is 
a solution of ${\cal L}_H \, U =0$ in $(t_1,
t_2) \times S^{2}$, for some $t_1 \leq  -  
\hat t_0$ and $t_2 \geq \hat t_0$, , (with $U =0$ at the boundary if $t_1>- \infty$ or $t_2 < +
\infty$) satisfying 
\[
| U |\leq (\cosh t)^{-\nu},
\]
for some $\nu > - \frac{1}{2}$, and with the property that, for every $t \in (t_1, t_2)$,  
$U(t, \cdot)$ is orthogonal to ${\cal V}^0$  in the $L^2$-sense on $S^{2}$, then $U \equiv 0$.
\label{pr:4-bis}
\end{proposition}
\noindent
{\bf Proof :} First we assume that  $U(t, \cdot)$ is not only orthogonal to ${\cal V}^0$  
in the $L^2$-sense on $S^{2}$ but also orthogonal to ${\cal V}^1$. Then (\ref{eq:est-5}) can
be improved into 
\[
D \geq 6 \, A.
\]
where $6$ corresponds to the next eigenvalue. And we can proceed as in Step 2 and Step 3 
 of the proof of the previous result, to show that $f=0$ and $v=0$.  Observe that, at this
point, there is no additional restriction needed.

\medskip

Therefore is just remains to prove the result when $U(t, \cdot)$ belongs to ${\cal V}^1$. 
In this case we are reduced to study some coupled system of ordinary differential equation.
Indeed,  we can decompose $f$ and $v$ over eigenfunctions  $\Delta^0$ and $\Delta^1$
associated to the eigenvalue $2$. Hence we can reduce to the case where 
\[
f= a \, \phi \qquad \mbox{and}\qquad v_{ex} = b \, d^0 \phi
\]
where $a$ and $b$ are scalar functions only depending on $t$ and $\phi$ is an 
eigenfunction of $\Delta^0$ corresponding to the eigenvalue $2$.. We obtain the following
system of ordinary differential equation
\begin{equation}
\left\{
\begin{array}{lllll}
\ds \ddot{a} - 2 \, a - \frac{9}{4}\, a - 4 \, \tanh (3t) \,  b  +   \frac{27}{4} \,  
\frac{a}{\cosh^2 (3t)} & = & 0 \\[3mm] 
\ds \ddot{b} - 2 \, b - \frac{1}{4}\, b - 2 \, \tanh (3t) \, a  +  \frac{7}{4} \,  
\frac{b}{\cosh^2 (3t)} & = & 0 .
\end{array}
\right.
\label{eq:syst}
\end{equation}
For the time being let us assume that $t_1=-\infty$ and $t_2= +\infty$ and show that 
any solution of (\ref{eq:syst}) which is bounded by $(\cosh s)^\nu$ vanishes. As already
mentioned, the asymptotic behavior of $a$ and $b$  at both $\pm \infty$ is governed by the
following sets of indicial roots
\[
\mu^{\pm}_1=\pm  \frac{5}{2}   \qquad \mbox{and}\qquad \nu^{\pm}_k=\pm  \frac{1}{2} 
\]
Observe that we know explicitely some solutions of (\ref{eq:syst}), namely the solutions 
corresponding to the Jacobi field $\Phi_t (\alpha,a)$. In particular, when $\alpha =0$, we
obtain the solution 
\[
\begin{array}{rlll}
a_1 & = &  (\sin (3s))^{- \frac{1}{6}} \, (\sin (3s) \, \cos  s - \sin s \, \cos (3s)) 
\\[3mm]
b_1 & = &   -  (\sin (3s))^{- \frac{1}{6}} \, \sin s 
\end{array}
\] 
Recall that $s$ is a function of $t$. Since $\sin (3s) = (\cosh (3t))^{-1}$ and $\cos (3s) 
=  - \tanh (3t)$, we can easily obtain the asymptotic behavior of this explicite solution.
Near $-\infty$ it is given by 
\[
\begin{array}{rlll}
a_1 & \sim & (\cosh t )^{- \frac{5}{2}}\\[3mm]
b_1 & \sim & - \frac{1}{3}\, (\cosh t)^{- \frac{5}{2}}
\end{array}
\]
and near $+\infty$ it is given by
\[
\begin{array}{rlll}
a_1 & \sim &  \sin \frac{\pi}{3} \, (\cosh t )^{- \frac{1}{2}}\\[3mm]
b_1 & \sim & - \sin \frac{\pi}{3} \, (\cosh t)^{- \frac{1}{2}}
\end{array}
\]
Now, assume  that we have a solution of (\ref{eq:syst}) which is bounded by $(\cosh t)^\nu$ 
for some $\nu < - \frac{1}{2}$. The inspection of the indicial roots shows that this solution
is bounded by a constant times  $(\cosh t)^{-\frac{5}{2}}$. However all solutions of
(\ref{eq:syst}) which are bounded by $(\cosh t)^{-\frac{5}{2}}$ have to be a multiple of the
solution $(a_1, b_1)$ described above.  Clearly these are not bounded by $(\cosh
t)^{-\frac{5}{2}}$ at $+\infty$ unless it is identically $0$.

\medskip

In order to complete the proof of the Proposition, we argue by contradiction and assume 
that the result is not true. There would exist sequences $(t'_i)_i \in [-\infty, 0]$ and
$(t_i'')_i \in [0, +\infty]$ tending to $-\infty$ and $+\infty$ respectively, and for each
$i$ a solution $(a_i, b_i)$ of (\ref{eq:syst}) defined in $(t'_i, t''_i)$ and bounded by a
$(\cosh s)^\nu$. These solutions have $0$ boundary data whenever $t_i'$ or $t_i''$ are
finite. 

\medskip

Our problem being linear we can assume that the solution $(a_i, b_i)$ is normalized in 
such a way that
\[
\sup_{(t_i', t_i'')} (\cosh t)^{-\nu} \, (|a_i|+-b_i|) =1
\]
Let $ t_i\in (t_i', t_i'')$ a point where the above maximum is achieved.  To begin with, 
observe that the sequence $t_i -t'_i$ remains bounded away from $0$. This is obvious if $t'_i
=-\infty$ and if not, this follows from the fact that since $(a_i, b_i)$  and therefore
$(\ddot a_i, \ddot b_i)$ are bounded by a constant (independent of $i$) times $(\cosh
t'_i)^\nu$ in $[t'_i, t'_i + 1]$ and since  $a_i=b_i =0$ at $t'_i$, standard ordinary
differential arguments show that $(\dot a_i, \dot b_i)$ is  also bounded  by a constant
(independent of $i$) times $(\cosh t'_i)^\nu$ in $[t'_i, t'_i +\frac{1}{2}]$. As a
consequence the above supremum cannot be achieved at a point which is too close to $t'_i$.
Similarly one proves that the sequence $t''_i-t_i$ also remains bounded away from $0$.

\medskip

\noindent
We now define the sequence of rescaled functions
\[
(\tilde a_i, \tilde b_i)  (t ) : =  (\cosh t_i)^\nu  \, (a_i, b_i) (t+ t_i).
\]

\medskip

\noindent
{\em Case 1 :} Assume that the sequence $t_i$ converges to $t_\infty \in {\R}$. 
After the extraction of some subsequences, if this is necessary, we may assume that the
sequence $(\tilde a_i, \tilde b_i) (\cdot - t_\infty)$ converges to some nontrivial solution
of (\ref{eq:syst}).  Furthermore this solution is bounded by a constant times $(\cosh
t)^\nu$. However, we have just proved that this is not possible. 

 \medskip

\noindent
{\em Case 2} : Assume that the sequence $t_i$ converges to $-\infty$. After the 
extraction of some subsequences, if this is necessary, we may assume that the sequence
$(\tilde a_i, \tilde b_i)$ converges to $(a_\infty, b_\infty)$ some nontrivial solution of  
\begin{equation}
\left\{
\begin{array}{lllll}
\ds \ddot{a} - \frac{17}{4}\, a + 4  \,  b   & = & 0 \\[3mm] 
\ds \ddot{b} -  \frac{9}{4}\, b + 2  \, a   & = & 0  
\end{array}
\right.
\label{eq:lklk}
\end{equation}
in some interval $(t_*,+\infty)$, with boundary condition $a_\infty =b_\infty= 0$, if  
$t_* : = \lim_{i\rightarrow \infty} t'_i-t_i$ is finite. Furthermore this solution is bounded
by a constant times $(\cosh t)^\nu$.   It is a simple exercice to show that (\ref{eq:lklk})
has no such solutions. 

\medskip

\noindent
{\em Case 3} : Finally, we assume that the sequence $t_i$ converges to $+\infty$. 
This case being similar to Case 2, we shall omit it.

\medskip

Since we have ruled out every possible case, the proof of the result is complete. 
\hfill $\Box$ 

\medskip

We will also need the following simpler result for the differential operator which 
appears in ${\cal L}_H$ when $t$ tends to $- \infty$.  Here no restriction are needed.
\begin{proposition}
Assume that $U$ is a solution of ${\cal L}_0 U =0$  in $(t_1, t_2) \times S^{n-1}$, 
(with $ U=0$ at the boundary if $t_1>- \infty$ or $t_2 < + \infty$). Assume that 
\[
|U|\leq (\cosh t)^{-\nu},
\]
for some $\nu > \frac{2-n}{2}$. Then $U \equiv 0$.
\label{pr:5}
\end{proposition}

\noindent
{\bf Proof :} This time the proof can be obtain be first decomposing $U=(f,v)$ over 
eigenspaces of $\Delta^0$ and $\Delta^1$ respectively and then compute explicitely the
solutions of the ordinary differential equation and finally show that $U=0$. 

\medskip

However, one can also proceed as in the former proof.  To begin with, let us assume 
that $n \geq 4$. We decompose $v = v_{ex}+ v_{coex}$ where $v_{ex}$ is an exact $1$-form on
$S^{n-1}$ and where $v_{coex}$ is a coexact $1$-form on $S^{n-1}$. 

\medskip

We multiply the equation ${\cal L}_H U$ by $(0, v_{coex})$ and integrate over $(t_1, t_2) 
\times S^{n-1}$. We obtain
\[
0  = \ds  \int  |\del_t v_{coex}|^2 + \int |d^1 v_{coex}|^2  +  \ds \frac{n^2}{4} \, \int  
|v_{coex}|^2 
\]
which already implies that $v_{coex}\equiv 0$. It therefore remains to prove that 
$(f, v_{ex})\equiv 0$.

\medskip

The proof is now quite involved and, in order to simplify the notations, we set
\[
\begin{array}{lllllllllll}
A :=  \ds \int  f^2,   & B  := \ds \int  |\del_t f|^2 , &  C  := \ds \int  |d^0 \,  f|^2  
\\[3mm]
A' :=  \ds \int  |v_{ex}|^2 ,   &  B'  :=  \ds \int |\del_t v_{ex}|^2, &  C'  :=  \ds \int 
|\delta^1 v_{ex}|^2.
\end{array}
\]

\medskip

We multiply the equation ${\cal L}_H U$ by $(f, v_{ex})$ and integrate. This time, 
we obtain
\[
\begin{array}{llll}
 \ds  \int ( | \del_t f |^2 +  |\del_t v_{ex}|^2) + \int ( |d^0 f|^2 + |\delta^1 v_{ex}|^2)+ 
 \ds \int \left(  \frac{n^2}{4} \,  f^2 + \frac{(n-4)^2}{4} \, |v_{ex}|^2\right) \\[3mm]
 =  -  2 \,\ds \int  (d^0\, f \wedge \star \, v_{ex} + f \, \delta^1 \,  v_{ex})
\end{array}
\]
Using with Cauchy-Schwarz inequality, we estimate
\begin{equation}
\begin{array}{rllll}
\ds  \left| \ds \int \, d^0 f \wedge \star \, v_{ex}  \right|   & \leq & \ds 
\left( \int |d^0 f|^2 \right)^{\frac{1}{2}}  \left(\int |v_{ex}|^2 \right)^{\frac{1}{2}}
\\[3mm]
\ds  \left| \ds \int \, f \, \delta^1 \, v_{ex} \right|  & \leq & \ds \left( \int 
|\delta^1 v_{ex}|^2 \right)^{\frac{1}{2}}  \left(\int f^2 \right)^{\frac{1}{2}}
\end{array}
\label{eq:est-5}
\end{equation}
To estimate the first RHS, we will use $2\, a \, b \leq  a^2 +   b^2$ and in order 
to estimate the second one, we use $2 \, a \, b \leq \frac{1}{2}\, a^2 + 2 \, b^2$.

\medskip

Collecting these, together with (\ref{eq:est-5}), we conclude that
\[
\ds  C + \frac{n^2}{4} A +  C' + \frac{(n-4)^2}{4} A'  \leq   \frac{1}{2} \, C' +
2 \, A + C +A'.
\]
To finish, we use the fact that, as in (\ref{eq:est-3}) we have 
\[
\int  |\delta^1 \, v_{ex} |^2 \geq (n-1) \, \int \, |v_{ex}|^2
\]
Otherwize stated that $C' \geq (n-1) \, A'$  to conclude that
\[
\ds  \frac{n^2-8}{4} A +  \frac{(n-4)^2+2n-6}{4} A'  \leq   0
\]
which proves the desired claim.  \hfill $\Box$ 

\section{Mapping properties of ${\cal L}_H$ on a half hyperbola}

As in \cite{Maz-Pac}, the analysis of the mapping properties of ${\cal L}_H$  is easy 
to do in some weighted H\"older spaces we are now going to define.
\begin{definition}
For all $\delta \in {\R}$ and for all $t_0 \in {\R}$, the space ${\cal C}^{k, 
\alpha}_{\delta}([t_0, +\infty) \times S^{n-1}; \Omega^0 \times \Omega^1)$ is defined to be
the space of  $U \in {\cal C}^{k, \alpha}( [t_0 , +\infty)\times S^{n-1};  \Omega^0 \times
\Omega^1)$ for which the following norm is finite
\[
\|U\|_{k, \alpha, \delta}:= \sup_{t \geq t_0} \,  |e^{-\delta s} \, U |_{k, \alpha \, 
([t, t+1]\times S^{n-1})}.
\]
Here $| \, \, \, |_{k, \alpha \, ([t, t+1]\times S^{n-1})}$ denotes  the usual H\"older 
norm in $[t, t+1]\times S^{n-1}$.
\end{definition}
Observe that $U=(f,v)$ where $v$ is a $1$-form, hence, in the above defined norm it is 
the coefficients of $v$ and the function $f$ which are estimated.

\medskip

To begin with, we investigate the mapping properties of ${\cal L}_H$ when defined
 between the above weighted spaces. These mapping properties crucially depend on the choice
of $\delta$. We prove the
\begin{proposition}
Assume that $\delta \in (\frac{2-n}{2}, \frac{n-2}{2})$ and $\alpha \in (0,1)$ are 
fixed. There exists some constant $c >0$ and, for all $t_0 \in {\R}$ (when $n=3$, $t_0$ has
to be chosen larger than $\hat t_0$ defined in Proposition~\ref{pr:4-bis}) , there exists an
operator 
\[
{\cal G}_{t_0} : {\cal C}^{0, \alpha}_\delta ([t_0, +\infty) \times S^{n-1};   
\Omega^0 \times \Omega^1) \longrightarrow  {\cal C}^{2, \alpha}_\delta ([t_0, +\infty) \times
S^{n-1};  \Omega^0 \times \Omega^1),
\]
such that, for all $V \in {\cal C}^{0, \alpha}_{\delta}([t_0, +\infty) \times S^{n-1};  
\Omega^0 \times \Omega^1)$, if for all $t >t_0$, $V(t, \cdot)$ is orthogonal to ${\cal V}^0$
in the $L^2$ sense on $S^{n-1}$, then  $U ={\cal G}_{t_0} V$ is the unique solution of 
\[
\left\{
\begin{array}{rlll}
{\cal L}_H U & = &  V \qquad  &\mbox{in} \qquad  [t_0, +\infty) \times S^{n-1}\\[3mm]
                  V & \in & {\cal V}^0  \qquad  & \mbox{on} \qquad  \{t_0\} \times S^{n-1} .
\end{array}
\right.
\]
Furthermore, $\|U\|_{2, \alpha, \delta} \leq c \, \|V\|_{0, \alpha, \delta} $.
\label{pr:6}
\end{proposition}

\noindent
{\bf Proof :} Uniqueness of ${\cal G}_{t_0}$ follows Proposition~\ref{pr:4}. We 
therefore concentrate our attention on the existence of ${\cal G}_{t_0}$ and the derivation
of the uniform estimate for the inverse. 

\medskip

Our problem being linear, we can assume that
\[
\sup_{(t', +\infty) \times S^{n-1}} |e^{- \delta t} \, V|\leq 1.
\]

Now, it follows from Proposition~\ref{pr:4} that, when restricted to the space of 
$U$ which are orthogonal to ${\cal V}^0$ in the $L^2$-sense on $S^{n-1}$, the operator ${\cal
L}_H$ is injective over $(t', t'')\times S^{n-1}$. As a consequence, for all $t''>t'+1$ we
are able to solve ${\cal L}_H U = V$, in $(t',t'')\times S^{n-1}$, with $ U =0$ on $\{t',
t''\}\times S^{n-1}$.

\medskip

We claim that, there exists some constant $c>0$ independent of $t''> t'+1$ and 
$t_0$ and of $V$ such that
\[
\sup_{(t',t'') \times S^{n-1}} |e^{- \delta t} \, U |\leq c.
\]
We argue by contradiction and assume that the result is not true. In this case, 
there would exist sequences $t_i'' > t_i'+1$, a sequence of functions $V_i$ satisfying
\[
\sup_{(t'_i, t''_i) \times S^{n-1}} |e^{ -\delta t} V_i | = 1,
\]
and a sequence $U_i$ of solutions of ${\cal L}_H (f_i, v_i) = (g_i, w_i)$, in 
$(t'_i,t''_i)\times S^{n-1}$, with $(f_i, v_i)=0$ on $\{t'_i, t''_i\}\times S^{n-1}$ such
that 
\[
A_i : = \sup_{(t'_i, t''_i) \times S^{n-1}} |e^{- \delta t} \, U_i | \longrightarrow +
\infty.
\]
Furthermore, $U_i (t,  \cdot)$ and $V_i (t, \cdot)$  are orthogonal in the $L^2$ sense to 
${\cal V}^0$ on $S^{n-1}$. Let us denote by $(t_i, \theta_i)\in (t'_i, t''_i)\times S^{n-1}$,
a point where the above supremum is achieved. We now distinguish a few cases according to the
behavior of the sequence $t_i$ (which, up to a subsequence can always be assumed to converge
in $[-\infty, +\infty]$). Up to some subsequence, we may also assume that the sequences
$t''_i -t_i$ (resp. $t_i-t'_i$) converges to $t^*\in (0, +\infty]$ (resp. to $t_* \in
[-\infty, 0)$). 

\noindent
Observe that the sequence $t_i-t'_i$ remains bounded away from $0$. Indeed,  since $U_i $ 
and   ${\cal L}_H \, U_i$ are bounded by a constant (independent of $i$) times $e^{\delta
t'_i} \, A_i$ in $[t'_i, t'_i + 1]\times S^{n-1}$ and since  $U_i =0$ on $\{t'_i \}\times
S^{n-1}$, standard elliptic estimates allow us to conclude that the partial derivative of
$U_i$ with respect to $t$ is also uniformly bounded by a constant times $e^{\delta t'_i} \,
A_i$ in $[t'_i, t'_i +\frac{1}{2}]\times S^{n-1}$. As a consequence the above supremum cannot
be achieved at a point which is too close to $t'_i$. Similarly one proves that the sequence
$t''_i-t_i$ also remains bounded away from $0$.

 \medskip

\noindent
We now define the sequence of rescaled functions
\[
\tilde{U}_i  (t, \theta) : = \frac{e^{-\delta t_i}}{A_i} \, U_i (t+ t_i, \theta).
\]

\medskip

\noindent
{\em Case 1 :} Assume that the sequence $t_i$ converges to $t_\infty \in {\R}$. After 
the extraction of some subsequences, if this is necessary, we may assume that the sequence
$\tilde{U}_i$ converges to some nontrivial solution of 
\[
{\cal L}_H \, U_\infty =0 ,
\]
in $(t_*,t^*)\times S^{n-1}$, with boundary condition $U_\infty=0$, if either $t_*$ or 
$t^*$ is finite. Furthermore
\begin{equation}
\sup_{(t_* ,  t^*)\times S^{n-1}}|e^{-\delta (t-t_\infty)} \, U_\infty | =1.
\label{eq:contr-1}
\end{equation}
and $U_\infty (t, \cdot)$ orthogonal in the $L^2$ sense to ${\cal V}^0$ on $S^{n-1}$.  
If $t_*=-\infty$, the inspection of the indicial roots shows that, any solution of the
homogeneous equation which is bounded by $e^{\delta t}$ for some $\delta \in (\frac{2-n}{2},
\frac{n-2}{2})$  is bounded by a constant times $e^{\frac{n-2}{2}t}$ at $-\infty$. Similarly,
if $t^*=+\infty$ we find that any solution of the homogeneous equation which is bounded by
$e^{\delta t}$ for some $\delta \in (\frac{2-n}{2}, \frac{n-2}{2})$  is bounded by a constant
times $e^{\frac{2-n}{2}t}$ at $-\infty$. But, applying Proposition~\ref{pr:4}, this implies
that $(f,v) =0$, contradicting (\ref{eq:contr-1}).

 \medskip

\noindent
{\em Case 2} : Assume that the sequence $t_i$ converges to $-\infty$. After the extraction 
of some subsequences, if this is necessary, we may assume that the sequence $\tilde{U}_i$
converges to some nontrivial solution of  
\begin{equation}
{\cal L}_0 U_\infty =0
\label{eq:equ-L0}
\end{equation}
in $(t_*,t^*)\times S^{n-1}$, with boundary condition $U_\infty = 0$, if either  $t_*$ or 
$t^*$ is finite. Furthermore
\begin{equation}
\sup_{(t_*, t^*) \times S^{n-1}} |e^{- \delta t} \, U_\infty | = 1,
\label{eq:contr-2}
\end{equation}
As already mentioned, the indicial roots corresponding to (\ref{eq:equ-L0}) are the same 
as the indicial roots corresponding to ${\cal L}_H$. Again, if $t_*=-\infty$, the inspection
of the indicial roots shows that, any solution of the homogeneous equation which is bounded
by $e^{\delta t}$ for some $\delta \in (\frac{2-n}{2}, \frac{n-2}{2})$  is bounded by a
constant times $e^{\frac{n-2}{2}t}$ at $-\infty$. Similarly, if $t^*=+\infty$ we find that
any solution of the homogeneous equation which is bounded by $e^{\delta t}$ for some $\delta
\in (\frac{2-n}{2}, \frac{n-2}{2})$  is bounded by a constant times $e^{\frac{2-n}{2}t}$ at
$-\infty$. But, applying Proposition~\ref{pr:5}, this implies that $U_\infty =0$,
contradicting (\ref{eq:contr-2}).

\medskip

\noindent
{\em Case 3} : Assume that the sequence $t_i$ converges to $+\infty$. This case being 
similar to Case 2, we shall omit it.

\medskip

Now that the proof of the claim is finished, we may pass to the limit $t''\rightarrow 
+\infty$ and obtain a solution of ${\cal L}_H U = V$, in $(t', +\infty)\times S^{n-1}$, with
$U=0$ on $\{t'\}\times S^{n-1}$, which satisfies
\[
\sup_{(t', +\infty)\times S^{n-1}} |e^{-\delta t}\, U |\leq c ,
\]
for some constant $c>0$ independent of $S$. To complete the proof of the Proposition, 
it suffices to apply Schauder's estimates in order to get the relevant estimates for all the
derivatives. \hfill $\Box$

\medskip

We now extend the right inverse ${\cal G}_{t_0}$ to the set of $V:= (g, 0)$ where $g$ 
only depends on $t$. For the sake of simplicity in the notations, we keep the same notation
for the right inverse since they are defined on orthogonal spaces. This is the content of the
following
\begin{proposition}
Assume that $\delta' < - \frac{n}{2}$ and $\alpha \in (0,1)$ are fixed. There exists some 
constant $c >0$ and, for all $t_0 \in {\R}$, there exists an operator 
\[
{\cal G}_{t_0} : {\cal C}^{0, \alpha}_{\delta'} ([t_0, +\infty) \times S^{n-1}; 
\Omega^0\times \Omega^1) \longrightarrow  {\cal C}^{2, \alpha}_{\delta'} ([t_0, +\infty)
\times S^{n-1}; \Omega^0\times \Omega^1),
\]
such that, for all $(g, 0) \in {\cal C}^{0, \alpha}_{\delta'}([t_0, +\infty) 
\times S^{n-1}; \Omega^0\times \Omega^1)$, if for all $t >t_0$, $g(t, \cdot)$ is constant on
$S^{n-1}$, then  $(f,0) ={\cal G}_{t_0} (g,0)$ is the unique solution of 
\begin{equation}
{\cal L}_H (f,0)  =  (g,0)
\label{eq:f0}
\end{equation}
in  $[t_0, +\infty) \times S^{n-1}$ which belongs to the space ${\cal C}^{2, 
\alpha}_{\delta'}([t_0, +\infty) \times S^{n-1}; \Omega^0\times \Omega^1)$.  Furthermore,
$\|(f,0)\|_{2, \alpha, \delta'} \leq c \, \|(g,0)\|_{0, \alpha, \delta'} $.
\label{pr:7}
\end{proposition}

\noindent
{\bf Proof :} The existence is easy. We already know an explicite solution of 
the homogeneous problem ${\cal L}_H (f_0, 0)=0$. This solution is given by 
\[
f_0 (t) : = (\cosh (nt))^{-\frac{1}{2}}.
\]
We define the solution of (\ref{eq:f0}) by
\[
f(t) : =  f_0(t) \, \int_t^{+\infty} (f_0 (\zeta))^{-2}\, \int_{\zeta}^{+\infty} 
f_0(\xi ) \, g(\xi)\, d\xi \, d\zeta .
\]
It is a simple exercise to show that this is a solution which is well defined  
and that the estimate is satisfied, since we have chosen $\delta' < -\frac{n}{2}$. \hfill
$\Box$

\medskip

We will also need the
\begin{proposition}
There exists $c>0$ such that, for all $t_0 \in {\R}$ and all $W \in {\cal C}^{2,
\alpha}(S^{n-1}; \Omega^0 \times \Omega^1)$,  there exists a unique $U_0 \in {\cal
C}^{2,\alpha}_{\frac{2-n}{2}}([t_0, +\infty) \times S^{n-1}; \Omega^0 \times \Omega^1)$
solution of 
\begin{equation}
\left\{ 
\begin{array}{rllll} 
{\cal L}_0 U_0 & =  & 0 \qquad & \mbox{\rm in}\quad (t_0 ,+\infty) \times S^{n-1}\\[3mm] 
                 U_0 & =  &  W \qquad & \mbox{\rm on} \quad \{ t_0 \}\times S^{n-1}.
\end{array}  
\right. 
\label{eq:poisson}
\end{equation}
Furthermore,  we have 
\[
|| U_0 ||_{2, \alpha, \frac{2-n}{2}}\leq c \, e^{\frac{n-2}{2} \, t_0}\, || W ||_{2,\alpha} ,
\]
\label{pr:8}
\end{proposition}
{\bf Proof :} Uniqueness of the solution follows Proposition~\ref{pr:5}. We therefore 
concentrate our attention on the existence of the solution and the derivation of the uniform
estimate for the inverse. 

\medskip

Our problem being linear, we can assume that
\[
\sup_{S^{n-1}} |e^{- \delta t} \, W |\leq 1.
\]

{\bf Step 1} Assume that the boundary data $W$ is orthogonal, in the $L^2$ sense, to
 ${\cal V}^1_{ex}$ the eigenspace of $(\Delta^0, \Delta^1)$ corresponding to the eigenvalue
$1-n$. We apply Proposition~\ref{pr:4} which implies that the operator ${\cal L}_0$ is
injective over $(t', t'')\times S^{n-1}$. As a consequence, for all $t''>t'+1$ we are able to
solve ${\cal L}_0 U = 0$, in $(t',t'')\times S^{n-1}$, with $U=W$ on $\{t'\}\times S^{n-1}$
and $U =0$ on $\{t''\}\times S^{n-1}$.

\medskip

We claim that, for any $\delta \in (- \frac{n}{2}, \frac{2-n}{2})$, there exists 
some constant $c>0$ independent of $t''> t'+1$ and $t_0$ and of $W$ such that
\[
\sup_{(t',t'') \times S^{n-1}} |e^{- \delta t} \, U |\leq c \, e^{- \delta t'}.
\]
We argue by contradiction and assume that the result is not true. In this case, 
there would exist sequences $t_i'' > t_i'+1$, a sequence of functions $W_i$ satisfying
\[
\sup_{S^{n-1}} |e^{ -\delta t} W_i | = 1,
\]
and a sequence $U_i$ of solutions of ${\cal L}_H U_i = 0$, in $(t'_i,t''_i)\times 
S^{n-1}$, with $U_i = W_i$ on $\{t'_i\}\times S^{n-1}$ and $U_i = 0$ on $\{t''_i\}\times
S^{n-1}$ such that 
\[
A_i : = \sup_{(t'_i, t''_i) \times S^{n-1}} |e^{\delta (t'_i -t)} \, U_i | 
\longrightarrow +\infty.
\]
Furthermore, $U_i(t,  \cdot)$  are orthogonal in the $L^2$ sense to ${\cal V}^1_{ex}$ on 
$S^{n-1}$. Let us denote by $(t_i, \theta_i)\in (t'_i, t''_i)\times S^{n-1}$, a point where
the above supremum is achieved. We now distinguish a few cases according to the behavior of
the sequence $t_i$ (which, up to a subsequence can always be assumed to converge in
$[-\infty, +\infty]$). Up to some subsequence, we may also assume that the sequences $t''_i
-t_i$ (resp. $t_i-t'_i$) converges to $t^*\in (0, +\infty]$ (resp. to $t_* \in [-\infty,
0)$). 

\noindent
As in the proof of Proposition~\ref{pr:6}, observe that the sequences $t_i-t'_i$ and  
$t''_i-t_i$ remain bounded away from $0$.

 \medskip

\noindent
We now define the sequence of rescaled functions
\[
\tilde{U}_i (t, \theta) : = \frac{e^{\delta (t'_i -t_i)}}{A_i} \, U_i(t+ t_i, \theta).
\]
Up to a subsequence, we can assume that the sequence $t''_i - t_i$ converges to $t^* 
\in (0, +\infty]$ and that $t_i -t'_i$ converges to $t_*\in [-\infty, 0)$. After the
extraction of some subsequences, if this is necessary, we may assume that the sequence
$\tilde U_i$ converges to some nontrivial solution of 
\begin{equation}
{\cal L}_0 \, U_\infty =0 ,
\end{equation}
in $(t_*,t^*)\times S^{n-1}$, with boundary condition $U_\infty =0$, if either $t_*$ or 
$t^*$ is finite. Furthermore
\begin{equation}
\sup_{(t_* ,  t^*)\times S^{n-1}}|e^{-\delta  t} \, U_\infty | =1.
\label{eq:contrad-1}
\end{equation}
and $U_\infty (t, \cdot) $ orthogonal in the $L^2$ sense to ${\cal V}_{ex}^1$ on $S^{n-1}$. 
If $t_*=-\infty$, the inspection of the indicial roots shows that, any solution of the
homogeneous equation which is bounded by $e^{\delta t}$ for some $\delta \in (- \frac{n}{2},
\frac{2-n}{2})$  is bounded by a constant times $e^{\frac{n}{2}t}$ at $-\infty$. Applying
Proposition~\ref{pr:6}, this implies that $U_\infty =0$, contradicting (\ref{eq:contrad-1}).

\medskip

Now that the proof of the claim is finished, we may pass to the limit $t''\rightarrow +
\infty$ and obtain a solution of ${\cal L}_0 U = 0$, in $(t', +\infty)\times S^{n-1}$, with
$U=W $ on $\{t'\}\times S^{n-1}$, which satisfies
\[
\sup_{(t', +\infty)\times S^{n-1}} |e^{-\delta t}\, U |\leq c ,
\]
for some constant $c>0$ independent of $t'$. Finally, Schauder's estimates yield the 
relevant estimates for all the derivatives. 

\noindent
{\bf Step 2}  Assume that the boundary data $W$ belongs to ${\cal V}^1_{ex}$ the 
eigenspace of $(\Delta^0, \Delta^1)$ corresponding to the eigenvalue $1-n$. In this case, the
partial differential equation ${\cal L}_0 U=0$ reduces to a finite number of coupled ordinary
differential equations of the form
\[
\left\{
\begin{array}{rllll}
\ddot{a} - (n-1) \, \dot{a} - \frac{n^2}{4} \, a + 2 (n-1) b & = & 0 \\[3mm]
\ddot{b} - (n-1) \, \dot{b} - \frac{(n-4)^2}{4} \, b + 2 a & = & 0
\end{array}
\right.
\]
with boundary data $a(t')=a_0$ and $b(t')=b_0$. Then, provided $\Delta^0 \phi = - 
(n-1) \, \phi$,  $U = (a \phi, b \, d^0 \phi)$ will be a solution of ${\cal L}_0 U =0$ in
$(t',+\infty) \times S^{n-1}$ and $U(t', \cdot)= (a_0 \, \phi, b_0\, d^0 \phi)$ on
$\{t'\}\times S^{n-1}$.

\medskip

But the above system can be solved explicitely
\[
\left\{
\begin{array}{rllll}
a & = & (n-1) \, A\, e^{-\frac{n+2}{2} (t-t')}+ B \, e^{\frac{2-n}{2} (t-t')} \\[3mm]
b & = & - \, A\, e^{-\frac{n+2}{2} (t-t')}+ B \, e^{\frac{2-n}{2} (t-t')} \\[3mm]
\end{array}
\right.
\]
and $(n-2) A=  a_0+b_0$, $(n-2) B = (n-1)b_0 -a_0$. It is easy to check that $\| U
 \|_{2, \alpha , \frac{3-n}{2}} \leq c$.  This completes the proof of the result. \hfill
$\Box$

\section{Structure of the mean curvature operator about the Hyperbola}

In order to understand the structure of the mean curvature operator for any surface 
close enough to the hyperbola, we go back to the variational definition of the minimal
surface. Let $V$ be a normal perturbation of the hyperbola parameterized by $X$. Recall that,
with our notations, 
\[
X =  \frac{e^{is}}{(\sin (ns))^{1/n}}\, \Theta
\]
and that the normal vector field can be taken to be 
\[
V = i \, e^{i(1-n)s}\, f\, \Theta + i \, e^{is}\, T
\]
We will denote by
\[
Y(t, \theta) := X(t, \theta) +V(t, \theta)
\]
First, let us compute the vectors which span the tangent space of the surface 
parameterized by $Y$.
\[
\del_t Y  = - ( \cosh (nt) )^{\frac{1}{n}} \, e^{i(1-n)s} \, \Theta +  e^{i(1-n)s}\, 
\left( \frac{(n-1)}{\cosh (nt)}  \, f + i \, \del_t f\right) \, \Theta - e^{is}\, \left(
\frac{1}{\cosh (nt)} \, T - i \del_t T \right)
\]
and, for all $j=1, \ldots, n-1$
\[
\del_j Y =  (\cosh (nt))^{\frac{1}{n}} \, e^{i\theta}\, \Theta_j + i \, e^{i(1-n)s}\, 
\del_{\theta_j} f\, \Theta - i \, e^{is}\, \left(\tanh (ns) + \frac{i}{\cosh (ns)}\right) \, 
\Theta_j + i \, e^{is}\, \del_{\theta_j} T
\]

The coefficients of the first fundamental form associated to $Y$ are given by
\[
\begin{array}{rlll}
|\del_t Y|^2  & =  &  (\cosh (nt))^{\frac{2}{n}} - 2 \, (n-1) \, (\cosh (nt))^{
\frac{1- n}{n}} \, f + |\del_t f|^2 + |\del_t T|^2 \\[3mm]
                    & +  &  (n-1)^2 \,  \cosh^2 (nt) \, f^2 + \cosh^2 (nt)  \,  |T|^2
\end{array}
\]
and, for all $j=1, \ldots, n-1$
\[
\begin{array}{rlll}
\del_t Y \cdot \del_j Y & =  & \left( - 2\, ( \cosh (nt))^{\frac{1-n}{n}} \, T 
\cdot \e_j + (\e_j f) \, \del_t f + \overline{\nabla}_{\e_j}T \cdot \del_t T  \right. \\[3mm]
		    &  +   &  \left. (n-2) \, \cosh^2 (nt) \, f \, T \cdot \e_j +   \tanh (nt) \, 
\del_t f \, T \cdot \e_j - \tanh (nt) \, f \, \del_t T \cdot \e_j \right) \, |\Theta_j|

\end{array}
\]

\[
\begin{array}{rllll}
|\del_j Y|^2 & = & \left( ( (\cosh (nt))^{\frac{2}{n}} + (\e_j f)^2 + |f|^2+ 
\overline{\nabla}_{\e_j} T|^2 +  2 \,  (\cosh (nt))^{\frac{1-n}{n}} \, f \right.\\[3mm]
	   &  -  & \left.  2 \, \tanh (nt) \, (\e_j f)  \, T \cdot \e_j  - 2 \, \tanh (nt) \, f \, 
(\overline{\nabla}_{\e_j} T \cdot \e_j)\right) \, |\Theta_j|^2
\end{array}
\]

Finally, if $j \neq k$ we have 
\[
\begin{array}{rllll}
\del_j Y \cdot \del_k Y & = & \left( (\e_j f) \, (\e_k f) + \overline{\nabla}_{\e_j} T 
\cdot \overline{\nabla}_{\e_k} T - \tanh (nt) \, ((\e_j f) \, T \cdot \e_k + (\e_k f) \, T
\cdot \e_j) \right. \\[3mm]
		    &  -  & \left.  \tanh (nt) \, f \,  (\overline{\nabla}_{\e_j}T \cdot \e_k +
\overline{\nabla}_{\e_k}T \cdot \e_j )\right) \, |\Theta_j|\, |\Theta_k|
\end{array}
\]

The exact value of the coefficients of the first fundamental form is not needed. 
We just observe that the first fundamental form of the surface parameterized by $Y$ is given
by
\[
{\mathbb I}_Y = {\mathbb I}_X + (\cosh (nt))^{\frac{1-n}{n}} \,L (f,T) +  Q (f,T)
\]
where ${\mathbb I}_X$ is the first fundamental form of the hyperbola in the variables 
$(t, \theta)$, namely
\[
{\mathbb I}_X  = (\cosh (nt))^{\frac{2}{n}} \, (dt^2 + (\cosh (nt))^{\frac{2}{n}} \, 
d\, \theta_i^2)
\]
And where $L$ is linear and $Q$ is quadratic in the variables $f, T, \del_t \, f, 
\del_t \, T$ and $\e_j f$, $\nabla^\tau_{\e_j} T$. Both $L$ and $Q$ have coefficients
depending on $t$ but they are all bounded and have bounded derivatives with respect to $t$.

\medskip

Once we have obtained the structure of the first fundamental form, it is a simple 
exercise to obtain the structure of the volume functional and then the structure of the
Euler-Lagrange equation. Hence, we conclude that the surface parameterized by $Y$ is minimal
if and only if $V$ is a solution of the following partial differential equation
\[
L_1 V +  (\cosh t)^{-1 - n}  \, \tilde{Q}_2 ( (\cosh t)^{-1}\, V) + (\cosh t)^{-1} 
\, \tilde{Q}_3 ( (\cosh t)^{-1}\, V) =0
\]
where $L_1$ is the linearized mean curvature operator, $\tilde{Q}_2$ is homogeneous of 
degree $2$ and where  $\tilde{Q}_3$ collects all the higher order terms. Observe that the
Taylor's coefficients of $\tilde{Q}_i$ are bounded functions of $t$ and so are the
derivatives of any order of these functions.

\medskip

When we conjugate  this operator, as we have done with $L_H$, we obtain the 
\begin{proposition}
The surface parameterized by 
\[
Y (t, \theta) = X(t, \theta)+ (\cosh (nt))^{\frac{2-n}{2n}}\, V(t, \theta)
\]
is minimal if and only if
\begin{equation}
{\cal L}_H \, V +  (\cosh t )^{-\frac{n}{2}}  \, Q_2 ((\cosh t)^{-\frac{n}{2}} \, V) + 
(\cosh t)^{\frac{n}{2}} \, Q_3 ((\cosh t)^{-\frac{n}{2}} \, V) =0
\end{equation}
where $Q_2$ is homogeneous of degree $2$ and where  $Q_3$ collects all the higher order 
terms. Observe that the Taylor's coefficients of $Q_i$ are bounded functions of $t$ and so
are the derivatives of any order of these functions.
\end{proposition}

\section{Minimal $n$-submanifolds close to the truncated hyperbola $H_I$}

For all $\e \in (0,1)$, $\rho_*>0$,  $\kappa >1$ and $ \beta \in [\frac{1}{\kappa}, 
\kappa ]$, we define $s_* \in (0, \frac{\pi}{n})$ and  $t_* \in {\R}$ by the identities
\[
\rho_* = ( n \, \beta \, \e)^{\frac{1}{n}} \, \frac{\cos s_*}{(\sin (ns_*))^{1/n}} ,
\]
and 
\[
e^{-n t_*} = \frac{\sin (ns_*)}{1 - \cos (ns_*)}.
\]
Hence, $t_* < 0$ for $\e$ small enough.

\medskip

Consider any normal perturbation of the rescaled hyperbola
\[
Y(t, \theta )= (n\, \beta \, \e)^{\frac{1}{n}} \left( X(t, \theta) + (\cosh (nt))^{
\frac{2-n}{2n}}V(t,\theta) \right) .
\]
As usual, we identify the normal vector field  $V$ with $U=(f,v)$. We have seen that 
$Y$ describes a minimal $n$-submanifold if and only if $U$ is a solution of 
\[
{\cal L}_H  U = (\cosh t)^{-\frac{n}{2}} \, Q_2 ((\cosh  t)^{-\frac{n}{2}} U) + (\cosh t)^{
\frac{n}{2}} \, Q_3 ((\cosh t)^{-\frac{n}{2}} U) .
\]

\medskip

We modify the normal bundle when $(t, \theta) \in [t_* , t_* + 2]\times S^{n-1}$ so 
hat it is now given by the "vertical plane" $\{i\, x\, : \, x\in {\R}^n\}$ for all $(t,
\theta) \in [t_* , t_* +1]\times S^{n-1}$. More precisely, instead of considering the normal
bundle spanned by 
\[
N_0 = i \, e^{i(1-n)s}\, \Theta \qquad \mbox{and}\qquad \forall j=1,\ldots, n 
\qquad N_j= i \, e^{i s}\, \frac{\Theta_j}{|\Theta_j|},
\]
we  want to choose the bundle spanned by the vectors
\[
\tilde N_0 : = i \,  \Theta \qquad \mbox{and}\qquad \forall j=1,\ldots, n \qquad 
\tilde N_j : = i \, \frac{\Theta_j}{|\Theta_j|},
\]
in $(t_*, t_*+\frac{1}{2})\times S^{n-1}$. As explained in \cite{Maz-Pac}, this 
modifies slightly the equation we have to solve into  
\begin{equation}
{\cal L}_H U  =  (\cosh t)^{-2n} \, L U+  (\cosh t)^{-\frac{n}{2}} \, \tilde Q_2 ((
\cosh (t)^{- \frac{n}{2}} U) +  (\cosh t)^{\frac{n}{2}} \, \tilde Q_3 ((\cosh
t)^{-\frac{n}{2}} U ) ,
\label{eq:nlpb}
\end{equation}
where the linear operator $L$ has coefficients which are bounded and supported in 
$[t_* , t_* + 2]\times S^{n-1}$ and where $\tilde Q_2$ and $\tilde Q_3$ enjoy the properties
of $Q_2$ and $Q_3$.

\medskip

\begin{definition}
We define $P_0$ to be the $L^2$ projection over the space orthogonal to ${\cal V}^0$  
in the $L^2$ sense on $S^{n-1}$.
\end{definition}
\medskip

We want to find $U$ solution of (\ref{eq:nlpb}) in $(t_*, +\infty) \times S^{n-1}$,
 boundary data $P_0 U = W$ on $\{t_* \} \times S^{n-1}$ where 
\[
\|W\|_{{\cal C}^{2, \alpha}} \leq \kappa \,  (\e \, \rho_*^{n})^{\frac{1}{2}},
\]
for some fixed constant $\kappa >0$. To begin with, let us solve
\begin{equation}
\left\{
\begin{array}{rlll}
{\cal L}_0 U_0 & =& 0  &\qquad \mbox{in}\qquad (t_*, +\infty) \times S^{n-1} \\[3mm]
                  U_0 & = & W &  \qquad \mbox{on}\qquad \{t_* \} \times S^{n-1}
\end{array}
\right.
\label{eq:lo=u}
\end{equation}
We already know from Proposition~\ref{pr:8} that there exists some constant $c>0$ such that
\[
\| U_0 \|_{2, \alpha, \frac{2-n}{2}} \leq c \, e^{\frac{n-2}{2}\, t_*}\, \| W\|_{2,\alpha}
\]

It remains to solve
\[
\begin{array}{rllll}
{\cal L}_H U & = &(\cosh t)^{-2n} \, L (U + U_0) -  {\cal L}_H U_0 +  (\cosh t)^{-
\frac{n}{2}} \, \tilde Q_2 ((\cosh (t)^{-\frac{n}{2}} (U + U_0))  \\[3mm]
                           & + & (\cosh t)^{\frac{n}{2}} \, \tilde Q_3 ((\cosh t)^{-
\frac{n}{2}} (U + U_0) 
\end{array}
\]
in $(t_*, +\infty) \times S^{n-1}$ with $P_0 \, U = 0$ on $\{t_* \} \times S^{n-1}$. 
We may then rewrite the previous equation as
\[
U = {\cal G}_{t_*} \,  {\cal N}(U),
\]
where by definition the nonlinear operator ${\cal N}$ is given by
\[
\begin{array}{rllll}
{\cal N} (U) & := & \ds   \left( (\cosh t)^{-2n} \, L (U + U_0)  - {\cal L}_H U_0+ 
 (\cosh t)^{-\frac{n}{2}} \, Q_2 ((\cosh (t)^{-\frac{n}{2}} (U + U_0)) \right. \\[3mm]
                              &  + & \left. \ds (\cosh t)^{\frac{n}{2}} \, Q_3 ((\cosh t)^{-
 \frac{n}{2}} (U + U_0) \right)
\end{array}
\]
and where $G_{t_*}$ is the right inverse constructed in Proposition~\ref{pr:6}. In order 
to solve this equation, we use a fixed point theorem for contraction mapping. This is the
content of the following
\begin{proposition}
Let $\delta \in (\frac{2-n}{2}, \frac{n}{6})$  be fixed. There exists $c_0 >0$ and for 
all $\kappa >0$, there exists $\e_0 >0$ such that, for all  $W$ which are orthogonal to
${\cal V}^0$ in the $L^2$ sense on $S^{n-1}$ and which satisfy 
\[
\| W\|_{2, \alpha}\leq \kappa \, (\e \, \rho_*^n)^{\frac{1}{2}},
\]
we have 
\[
\|{\cal G}_{t_*} \, {\cal N} (0)\|_{2, \alpha, \delta}\leq \frac{c_0}{2}\, (\cosh t_*)^{
\frac{2-n}{2}} \, \|W\|_{2, \alpha}
\]
and
\[
\|{\cal G}_{t_*} \,  ({\cal N} (U_2) - {\cal N} (U_1))\|_{2, \alpha, \delta}\leq 
\frac{1}{2}\, \|U_2 -U_1\|_{2, \alpha, \delta}
\]
for all $U_1, U_2 \in {\cal C}^{2, \alpha}_\delta ([t_*, +\infty) \times S^{n-1})$ 
such that $\|U_i\|_{2, \alpha , \delta} \leq c_0 \, (\cosh t_*)^{\frac{2-n}{2}} \, \|W\|_{2,
\alpha}$.

\medskip

In particular ${\cal G}_{t_*} \, {\cal N}$ has a unique fixed point in the ball of 
radius $c_0 \, (\cosh t_*)^{\frac{2-n}{2}} \, \|W\|_{2, \alpha}$ in ${\cal C}^{2,
\alpha}_\delta ([t_*, +\infty) \times S^{n-1}; \Omega^0 \times \Omega^1)$.
\label{pr:10}
\end{proposition}
{\bf Proof:}  To begin with let us estimate
\[
\|  \left( (\cosh t)^{-2n} \, L  U_0  \right)\|_{2, \alpha, \delta} \leq c \, (
\cosh t_*)^{\delta - 2n} \, \| W\|_{2, \alpha}
\]
Now
\[
{\cal L}_H U_0  = ({\cal L}_H  - {\cal L}_0 ) U_0
\]
hence 
\[
\|   {\cal L}_H U_0 \|_{2, \alpha}  \leq c \, \left( (\cosh t_*)^{\delta-n} + (
\cosh t_*)^{\frac{2-n}{2}} \right)  \, \|W\|_{2, \alpha}.
\]
Now, using the properties of $\tilde{Q}_2$ and $\tilde{Q}_3$, we estimate for 
$\e$ small enough
\[
\| (\cosh t)^{-\frac{n}{2}} \, \tilde Q_2 ((\cosh t)^{-\frac{n}{2}} U_0) \|_{2, 
\alpha, \delta} \leq c \, \left( (\cosh t_*)^{\delta-\frac{3n}{2}} + (\cosh t_*)^{2-n}
\right)  \, \|W\|^2_{2, \alpha}
\]
and 
\[
\|   (\cosh t)^{\frac{n}{2}} \, \tilde Q_3 ((\cosh  t)^{-\frac{n}{2}} U_0)  \|_{2, 
\alpha, \delta} \leq c \, \left( (\cosh t_*)^{\delta- n} + (\cosh t_*)^{\frac{3(2-n)}{2}}
\right)  \, \|W\|^3_{2, \alpha}
\]
We can now use the result of Proposition~\ref{pr:6} to conclude that
\begin{equation}
\|{\cal G}_{t_*} \, P_0 \,  {\cal N} (0)\|_{2, \alpha, \delta}\leq c \, (\cosh t_*)^{
\frac{2-n}{2}} \, \|W\|_{2, \alpha}
\label{eq:f-estim}
\end{equation}
where we recall that $P_0$  to be the $L^2$ projection over the space orthogonal to 
${\cal V}^0$  in the $L^2$ sense on $S^{n-1}$.

\medskip

Using the explicit representation of ${\cal G}_{t_*}$ when acting on functions only 
depending on $t$, we obtain
\[
\| {\cal  G}_{t_*} \, (I- P_0) \, \left( (\cosh t)^{-2n} \, L   U_0  \right)\|_{2, 
\alpha, \delta} \leq c \, (\cosh t_*)^{\delta - \frac{3n}{2}} \, \| W \|_{2, \alpha}
\]
Now
\[
(I- P_0) \, {\cal L}_H U_0  = 0
\]
Finally, we estimate for all $\e$ small enough
\[
\begin{array}{rlll}
\| {\cal G}_{t_*} \, (I- P_0) \,  \left(  (\cosh t)^{-\frac{n}{2}} \, \tilde Q_2 
((\cosh t)^{-\frac{n}{2}} U_0) \right) \|_{2, \alpha, \delta} & \leq & c \, \left( (\cosh
t_*)^{\delta+ \frac{4-n}{2}}  \right. \\[3mm]
											& + &  \left. (\cosh t_*)^{2-n} \right)  \, \| W \|^2_{2, \alpha}
\end{array}
\]
and 
\[
\begin{array}{rllll}
\| {\cal G}_{t_*} \, (I- P_0) \,  \left(  (\cosh t)^{\frac{n}{2}} \, \tilde Q_3 
((\cosh t)^{-\frac{n}{2}} U_0) \right) \|_{2, \alpha, \delta}&  \leq & c \, \left( (\cosh
t_*)^{\delta+ 3- n} \right. \\[3mm]
							& + & \left.  (\cosh t_*)^{\frac{3(2-n)}{2}} \right)  \, \| W \|^3_{2, \alpha}
\end{array}
\]

Collecting these estimates together with (\ref{eq:f-estim}), we obtain 
\[
\|{\cal N} (0)\|_{2, \alpha, \delta}\leq c \, (\cosh t_*)^{\frac{2-n}{2}} \, \| W 
\|_{2, \alpha}
\]
provided $\e$ is chosen small enough. We fix $c_0 >0$ to be equal to twice the 
constant which appears on the right hand side of this estimate. The second estimate requires
that $\delta \in (\frac{2-n}{2}, \frac{n}{6})$. \hfill $\Box$

\section{Family of half hyperbola parameterized by their boundary}

We summarize what we have obtained so far. Given $\kappa >1$, and $\rho_* >0$, 
there exists $\e_0 >0$ such that for all $\beta \in [\frac{1}{\kappa}, \kappa]$, for all $\e
\in (0, \e_0)$ and all $\Phi \in {\cal C}^{2, \alpha}(S^{n-1}; {\R}^n)$ which is orthogonal
to $\Theta$ in the $L^2$ sense and which satisfies
\[
\|\Phi\|_{2, \alpha} \leq \kappa \, \rho_*^2 \, \e,
\]
we have obtained a minimal $n$-submanifold which is ${\cal C}^{2, \alpha}_\delta$ 
close to the truncated hyperbola, more precisely this submanifold is parameterized by
\[
Y  : (t,\theta) \in [t_*, +\infty) \times S^{n-1} \longrightarrow  (n\, \beta \, 
\e)^{\frac{1}{n}} \left( X(t, \theta) + (\cosh (nt))^{\frac{2-n}{2n}}V(t,\theta) \right)
\]
where $V =U+U_0 \in {\cal C}^{2, \alpha}_\delta ([t_*,+\infty) \times S^{n-1})$. 
Here $U$ is the solution of (\ref{eq:lo=u}) and $U_0$ is the solution of (\ref{eq:lo=u})
obtained in the previous Proposition with $W$ given by
\[
W : = (n\, \beta \, \e)^{- \frac{1}{n}} \, (\cosh (nt_*))^{\frac{n-2}{2n}} \, \Phi .
\]

\medskip

Now, we would like to describe the boundary of this $n$-submanifold. To this aim, 
let us define 
\[
r := (n \beta \e)^{\frac{1}{n}} \frac{\cos s}{( \sin (ns))^{\frac{1}{n}}}.
\]
Thanks to (\ref{eq:expansion}), we see that the above minimal $n$-submanifold can 
be parameterized, in some neighborhood of its boundary as
\[
(r, \theta) \longrightarrow  r \, \Theta + i \left( \e \,  \beta \,  r^{1-n} \, 
\Theta + F_0 +F )\right) (1+ {\cal O}(\e^2 r^{-2n}))
\]
Where $F_0$ is the harmonic extension of the boundary data $\Phi$ in $B_{\rho_*}$ 
and where 
\[
\| F(\rho_*, \cdot) \|_{{\cal C}^{2, \alpha} (S^{n-1})} + \| \rho_* \, \del_r F(
\rho_*, \cdot) \|_{{\cal C}^{1, \alpha} (S^{n-1})}\leq c \,  (\cosh
t_*)^{\frac{2-n}{2}-\delta} \, \|\Phi\|_{2, \alpha}.
\]
In other words we can also write this $n$-submanifold near its boundary as the graph of 
\[
 (r, \theta)  \longrightarrow   \e \, \beta  \, r^{1-n} \, \Theta + F_0 + \tilde{F},
\]
where $\tilde{F}$ depends smoothly on $\beta$ and $\Phi$ and satisfies
\[
\| \tilde F (\rho_*, \cdot) \|_{{\cal C}^{2, \alpha} (S^{n-1})} + \| \rho_* \, 
\del_r \tilde F (\rho_*, \cdot) \|_{{\cal C}^{1, \alpha} (S^{n-1})}\leq c_\kappa  \, \left(
\e^3 \rho_*^{1-3n} + \e \, \rho_* \, (\e^{-\frac{1}{n}} \, 
\rho_*)^{\frac{2-n}{2}-\delta}\right) ,
\]
where the constant $c_\kappa >0$ only depends on $\kappa$.

\medskip

For further use it will be convenient to define 
\begin{definition}
For all $\Phi \in {\cal C}^{2, \alpha} (S^{n-1}; {\R}^n)$, we define ${\cal P}_{
int}(\Phi)$ to be equal to $\del_r F$ where $F$ is the harmonic extension of $\Phi$ in the
unit ball of ${\R}^n$.
\end{definition}

Finally, given any ${\cal R}\in {\cal O}(n)$ we consider the image of the above 
defined family of minimal $n$-submanifolds by
\[
x+iy \in {\C}^n \longrightarrow x + i \, {\cal R} \, y \in {\C}^n.
\]
We denote by $H_{\e, {\cal R}, \rho_*} (\beta, \Phi )$ this submanifold.

\section{Mean curvature operator for normal graphs over ${\R}^n$}

Let a normal graph over the $x$ space in ${\C}^n$, namely $\Omega \subset {\R}^n \ni
 x \longrightarrow  x + i \, F(x)\in {\C}^n$. The first fundamental form is given by
\[
{\mathbb I} = \sum_j dx_j \, dx_j + 2 \, \sum_{j< j'} \del_{x_j} F \cdot \del_{x_j'} 
F \,  dx_j \, dx_{j'} 
\]

It is an easy exercise to derive the equation which ensures that the graph of $F$ will 
be minimal. The exact expression of this equation will not be needed but we only need its
structure 
\[
\Delta  F + \mbox{div} \,  \left( Q_3 (\nabla F)\right) =0
\]
where $Q_3$ is analytical and has coefficients which do not depend on $x$. Furthermore, 
$Q_3$, $\nabla Q_3$ and $\nabla^2 Q_3$ all vanish at $0$.

\section{Mapping properties of the Laplace operator in ${\R}^n$}

Choose $\rho_0 >0$ such that, for all $j \neq j'$, we have $B(x_j, \rho_0) \cap B(x_{j'}, 
\rho_0) = \emptyset$ and also $B(x_j, \rho_0) \subset B(0, \rho_0^{-1})$. 

\medskip

The weighted space we will be working with is defined by
\begin{definition}
Given $\nu, \mu \in {\R}$ and $\alpha \in (0,1)$. The space ${\cal C}^{k, \alpha}_{\nu, \mu} 
({\R}^n \setminus \{x_1, \ldots, x_N\}; {\R}^n)$ is defined to be the space of functions $F
\in {\cal C}^{k, \alpha}_{cal} ({\R}^n \setminus \{x_1, \ldots, x_N\}; {\R}^n)$ for which,
the following norm is finite
\[
\| F\|_{k, \alpha, \nu, \mu} : =  \sum_j \sup_{r\in (0, \rho_0)} r^{-\nu}\,  [F]_{k, \alpha, 
B_{r}\setminus B_{r/2}(x_j)} + |F|_{k, \alpha, B_{2/\rho_0}\setminus \cup_j
B_{\rho_0/2}(x_j)} + \sup_{r \geq 1/\rho_0} r^{-\mu} \, [F]_{k, \alpha, B_{2r}\setminus B_r},
\]
where 
\[
[F]_{k, \alpha, B_{r}\setminus B_{r/2}} : = \sum_{j=0}^k r^{j}\, \sup_{B_{r}\setminus 
B_{r/2}} |\nabla^j F| + r^{k+\alpha} \, \sup_{x,y \in B_{r}\setminus B_{r/2}}\frac{|\nabla^k
F(x)- \nabla^k F (y)|}{|x-y|^\alpha}.
\]
and where $|\, \, |_{k, \alpha, \Omega}$ is the usual H\"older norm in $\Omega$.
\end{definition}

\begin{proposition}
Assume that $2-n < \nu < 0$ and also that $1-n < \mu < 2-n$. Then
\[
\Delta : {\cal C}^{2, \alpha}_{\nu, \mu}({\R}^n \setminus \{x_1, \ldots, x_N\}; {\R}^n) 
\oplus (1+ |x|^2)^{\frac{2-n}{2}}\times {\R}^n \longrightarrow {\cal C}^{0, \alpha}_{\nu-2,
\mu-2}({\R}^n \setminus \{x_1, \ldots, x_N\}; {\R}^n) ,
\]
is an isomorphism.
\end{proposition}
 {\bf Proof :} For all $H \in {\cal C}^{0, \alpha}_{\nu-2, \mu-2}({\R}^n \setminus \{x_1, 
\ldots, x_N\}; {\R}^n)$, $H \in L^1 ({\R}^n)$. Hence there exists $F$ weak solution of
$\Delta F=H$ in ${\R}^n$. Now, the function $x \rightarrow |x-x_j|^\nu$ can be used as a
barrier function to prove that $F$ is bounded by a constant times $\|H\|_{0, \alpha, \nu-2,
\mu-2}$ times $|x-x_j|^\nu$ in each $B(x_j, \rho_0)$. 

\medskip

Moreover, it is classical to prove that $F$ is bounded by a constant times $\|H\|_{0, 
\alpha, \nu-2, \mu-2}$ times $|x-x_j|^{2-n}$ outside $\cup_j B(x_j, \rho_0)$. The estimates
for the derivatives follow from rescaled Schauder's estimates. \hfill $\Box$ 

\medskip

As above, we assume that $\rho_0 >0$ is chosen such that, for all $j \neq j'$, we have 
$B(x_j, \rho_0) \cap B(x_{j'}, \rho_0) = \emptyset$. Using the maximum principle and the
previous result, it is a simple exercise to show that
\begin{proposition}
Assume that $2-n < \nu < 0$ and also that $1-n < \mu < 2-n$. Then, for all $\rho \in (0, 
\rho_0)$
\[
\Delta : {\cal C}^{2, \alpha}_{\nu, \mu, {\cal D}}({\R}^n \setminus \cup_{j=1}^N B(x_j, 
\rho); {\R}^n ) \oplus (1+ |x|^2)^{\frac{2-n}{2}}\times {\R}^n \longrightarrow {\cal C}^{0,
\alpha}_{\nu-2, \mu-2}({\R}^n \setminus \cup_{j=1}^N B(x_j, \rho); {\R}^n),
\]
is an isomorphism. Moreover, there exists $c >0$ such that, for all $\rho \in (0, \rho_0)$ 
the norm of its inverse ${\cal G}_{\rho}$ is bounded by $c$.
\label{pr:12}
\end{proposition}
Here the subscript ${\cal D}$ refers to  the fact that all functions have $0$ boundary
data on each $\del B (x_j, \rho)$. 

\section{Expansion of Green's function}

Assume that $x_1, \ldots, x_N \in {\R}^n$,  $\alpha_1, \ldots, \alpha_N \in {\R}$, 
${\cal A}_0\in M_n ({\R})$ and ${\cal R}_1, \ldots, {\cal R}_N \in O(n)$ are given, we set 
\[
G(x)  = \sum_{j=1}^N \alpha_j {\cal R}_j \left( \frac{x-x_j}{|x-x_j|^n}\right)  + 
{\cal A}_0 \, x .
\]
It is an easy exercise to perform an expansion of $G$ near one the point $x_{j_0}$. 
We get
\[
\begin{array}{rllll}
G(x)    &  =  & \ds  \alpha_{j_0} \, {\cal R}_{j_0} \left( \frac{x-x_{j_0}}{|x-x_{j_0}|^n}
\right)  -  \ds \left( \sum_{j\neq j_0} \alpha_j \, {\cal R}_j \left(
\frac{x_j}{|x_j|^n}\right) - {\cal A}_0  \, x_{j_0} \right) + {\cal A}_0 \,  (x-x_{j_0})
\\[3mm]
           &  +  & \ds \sum_{j\neq j_0} \frac{\alpha_j}{|x_j-x_{j_0}|^{n}}  \, {\cal R}_j 
\left( x-x_{j_0} -  n \, ((x-x_{j_0}) \, \cdot \,  (x_j-x_{j_0})) \,
\frac{x_j-x_{j_0}}{|x_j-x_{j_0}|^2}\right) \\[3mm]
           &  +  &  \ds {\cal O} ( |x-x_{j_0}|^2).
\end{array}
\]
We set 
\[
 \rho \, \Theta : = x-x_{j_0}\qquad \mbox{and, for all $j \neq j'$} \qquad \xi_{jj'}:=  
\frac{x_j-x_{j'}}{|x_j -x_{j'}|},
\] 
and we define $\gamma_{jj}=0$ for all $j$ and 
\begin{equation}
\gamma_{jj'} :=  \frac{1}{|x_j-x_{j'}|^{n}} \left( {\int}_{S^{n-1}} {\cal R}_j  \Theta \cdot 
{\cal R}_{j'} \Theta \, d\theta -  n \,  \int_{S^{n-1}} (\Theta \cdot {\cal R}_j 
\xi_{j,j'})\, (\Theta \cdot {\cal R}_{j'} \xi_{j,j'}) \, d\theta \right),
\label{eq:matrix}
\end{equation}
for all $j\neq j'$, and 
\[
\lambda_j  =  - \int_{S^{n-1}} {\cal A}_0 \Theta \cdot {\cal R}_{j}\Theta \, d\theta
\]
so that we can write
\begin{equation}
\begin{array}{rlll}
G(x)   & = &  \ds \left( \alpha_{j_0} \, \rho^{1-n} +  \frac{1}{\omega_n} \, ( 
\sum_{j\neq j_0} \alpha_j  \, a_{j j_0} -  \lambda_{j_0} ) \, \rho \right)  \,  {\cal
R}_{j_0}\Theta \\[3mm]
          & -  & \ds \left( \sum_{j\neq j_0} \alpha_j \, {\cal R}_j \left( 
\frac{x_j}{|x_j|^n}\right) -{\cal A}_0 x_{j_0}\right) + {\cal O} (\rho^2) + {\cal O}_{\perp} 
(\rho)  
\end{array}
\label{eq:expansion-bis}
\end{equation}
where $ {\cal O}_{\perp}  (\rho) $ is orthogonal to ${\cal R}_{j_0} \Theta$ in 
the $L^2$ sense and where $\omega_n := |S^{n-1}|$.

\section{Minimal graphs}

Assume that $x_1, \ldots, x_N \in {\R}^n$,  $\alpha_1, \ldots, \alpha_N \in {\R}$, 
${\cal A}_0 \in M_n ({\R})$ and ${\cal R}_1, \ldots, {\cal R}_N \in O(n)$ are given, we have
defined 
\[
G(x) : = \sum_{j=1}^N \alpha_j {\cal R}_j \left( \frac{x-x_j}{|x-x_j|^n}\right)  + 
{\cal A}_0 \, x
\]
We would like to solve
\[
\left\{
\begin{array}{rllll}
\Delta F + \mbox{div} \left( Q_3 (\nabla F)\right) & = & 0 \qquad & \mbox{in}\qquad 
{\R}^n \setminus \cup_j B(x_j, \rho_*) \\[3mm]
                                                                            F & = & \e \, 
G + {\cal R}_j \, \Phi_j  \qquad & \mbox{on} \qquad \del B(x_j,\rho_*) 
\end{array}
\right.
\]
where $\Phi_j$ is small.

\medskip

To begin with, let us define $F_j$ to be the solution of 
\[
\left\{
\begin{array}{rllll}
\Delta F_j & = & 0 \qquad & \mbox{in}\qquad {\R}^n \setminus \cup_j B(x_j, \rho_*) \\[3mm]
           F_j & =  &  {\cal R}_j \, \Phi_j  \qquad & \mbox{on} \qquad \del B(x_j, \rho_*) 
\end{array}
\right.
\]
Using barrier arguments, we obtain the 
\begin{lemma}
There exists a constant $c >0$ such that
\[
\| F_j \|_{2, \alpha, 2-n, 2-n}\leq c \, \rho_*^{n-2} \, \| \Phi_j\|_{2, \alpha}.
\]
\label{eq:estimfj}
\end{lemma}
{\bf Proof :} Observe that 
\[
x \longrightarrow \left( \frac{|x-x_j|}{\rho_*}\right)^{2-n} \, \| \Phi_j\|_{2, \alpha}.
\]
is a supersolution for our problem, hence we already have
\[
| F_j | \leq  \left( \frac{|x-x_j|}{\rho_*}\right)^{2-n} \, \| \Phi_j\|_{2, \alpha}.
\]
The estimates for the other derivatives follow from Schauder's estimates. \hfill $\Box$

\medskip

As above, we choose $\rho_0 >0$ such that, for all $j \neq j'$, we have $B(x_j, \rho_0) 
\cap B(x_{j'}, \rho_0) = \emptyset$ and also $B(x_j, \rho_0) \subset B(0, \rho_0^{-1})$.
Using cutoff functions, we can restrict the definition of $F_j$ to $B(x_j,\rho_0)$ and sum
all these functions to obtain a mapping 
\[
\tilde{F}_0 : = \sum_{j} \chi_{\rho_0} (x-x_j) \, F_j (x)
\]
where $\chi_{\rho_0}$ is equal to $1$ in $B(0, \rho_0/2)$ and equal to $0$ outside $B(0, 
\rho_0)$.
 
\medskip

It is an easy task to evaluated the mean curvature of the graph of $F = \e \, G + \tilde 
F_0$. This is the contain of the following 
\begin{lemma}
Assume that $\nu \in (2-n, 0)$ and that $\mu \in (1-n, 2-n)$. There exists a constant $c 
>0$ such that
\[
\|\Delta F + \mbox{div} \left( Q_3 (\nabla F)\right)\|_{0, \alpha, \nu-2, \mu-2}\leq c \, 
\left( \e^3 \rho_*^{1-3n-\nu}+ \sup_j \|\Phi_j\|_{2, \alpha} \, \rho_*^{n-2}\right)
\]
\label{le:estimh}
\end{lemma}
{\bf Proof :} In order to evaluate the mean curvature of the graph of $F$, it suffices to
 evaluate $\Delta F + \mbox{div} \left( Q_3 (\nabla F)\right)$. We use the fact that $\Delta
F =0$ away from each $B(x_j, \rho_0)$ and also in $B(x_j, \frac{\rho_0}{2})\setminus B(x_j,
\rho_*)$ while $\Delta F$ and all its derivatives are bounded by a constant times $
\rho_*^{n-2}\,  \| \Phi_j\|_{2, \alpha}$ in each $B(x_j, \rho_0) \setminus B(x_j,
\frac{\rho_0}{2})$. Hence 
\[
\| \Delta F  \|_{2, \alpha, \nu-2 , \mu-2}\leq c \,  \rho_*^{n-2}\,  \sup_j \| \Phi_j\|_{2,
 \alpha}
\]
In order to estimate the norm of $\mbox{div} \left( Q_3 (\nabla F)\right)$ we first observe
 that the mean curvature of the graph of the function $F_0 : x \rightarrow \e \, {\cal A}_0
\, x$ is $0$ hence $\mbox{div} \left( Q_3 (\nabla F)\right) = \mbox{div} \left( Q_3 (\nabla
F)\right)-\mbox{div} \left( Q_3 (\nabla F_0)\right)$. We can thus evaluate 
\[
\|  \mbox{div} \left( Q_3 (\nabla F)\right) \|_{2, \alpha, \nu-2 , \mu-2}\leq c \, 
\left( \e^3 \rho_*^{1-3n-\nu}+ \sup_j \| \Phi_j\|_{2, \alpha} \, \rho_*^{n-2}\right)
\]
This ends the proof of the Lemma. \hfill $\Box$

\medskip

The graph of $F= \e \, G + \tilde F_0$ over the horizontal plane $\Pi_0$ can also be
 viewed as a normal graph over the submanifold $\tilde \Pi_\e$ which is the graph of $x
\rightarrow x+ i  \, (1- \chi_{\frac{2}{\rho_0}}) \, \e \, {\cal A}_0 \, x$, where 
$\chi_{\frac{2}{\rho_0}}$ is equal to $0$ in $B(0, \frac{1}{\rho_0})$ and equal to $0$
outside $B(0, \frac{2}{\rho_0})$. By construction, this submanifold $\tilde \Pi_0$ is equal
to $\Pi_0$ in $B(0, \frac{1}{\rho_0})$ and is equal to the graph of the mapping $F_0$ away
from $B(0, \frac{2}{\rho_0})$. Hence we can say that the vertical graph of $F = \e  \, G +
\tilde F_0$ over $\Pi_0$ is a normal graph of some mapping $\tilde F$ over $\tilde \Pi_0$. 

\medskip

Now we want to perturb this normal graph in order to obtain a minimal submanifold.  
The equation we have to solve now reads 
\[
\left\{
\begin{array}{rllll}
\Delta (F+ \tilde F) + \e \, L (F+\tilde F) +  \mbox{div} \left( \tilde Q_3 (\nabla 
(F+ \tilde F))\right) & = & 0 \qquad & \mbox{in}\qquad {\R}^n \setminus \cup_j B(x_j, \rho_*)
\\[3mm]
                                                                                                                                                          F & = & 0  \qquad & \mbox{on} \qquad \del B(x_j, \rho_*)  .
\end{array}
\right.
\]
Where $L$ is a second order linear differential operator which takes into account the fact 
that  $\tilde \Pi_0$ has some region  where it is not planar.  The coefficients of $L$ are
bounded independently of $L$ and can be assumed to be supported in $B(0, \frac{2}{\rho_*})
\setminus B(0, \frac{1}{\rho_*})$.

\medskip

Making use of Proposition~\ref{pr:12}, we can rewrite this equation as 
\[
F = - {\cal G}_{\rho_*} \,  \left( \Delta \tilde F + \e \, L (F+\tilde F) +  \mbox{div} 
\left( \tilde Q_3 (\nabla (F+ \tilde F))\right)  \right) .
\]
For the sake of simplicity, let us call ${\cal M}(F)$ this operator which is defined on the 
space 
\[
{\cal E}_{\rho_*} :=  \left[ {\cal C}^{2, \alpha}_{\nu, \mu} ({\R}^n \setminus \cup_j B 
(x_j , \rho_*); {\R}^n) \oplus  (1+|x|^2)^{\frac{2-n}{2}}\times {\R}^n\right]_{\cal D},
\]
which is naturally endowed with the product norm, with values into ${\cal C}^{2, \alpha}_
{\nu-2, \mu-2} ({\R}^n \setminus \cup_j B(x_j , \rho_*); {\R}^n) $.

\medskip

The existence of a fixed point for this operator is the content of the following 
\begin{proposition}
Let $\nu \in (2-n, 0)$, $\mu \in (1-n, 2-n)$ and $\kappa >0$ be fixed. There exists $c_0 >0$
 and  $\e_0 >0$ such that for all  $\e \in (0, \e_0)$ and for all $\Phi_j \in {\cal C}^{2,
\alpha}(S^{n-1}; {\R}^n)$  which satisfy $\| \Phi_j \|_{2, \alpha}\leq \kappa \, \e \,
\rho_*^2$, we have 
\[
\|{\cal M} (0)\|_{\cal E}\leq \frac{c_0}{2}\, ( \e^3 \, \rho_*^{1-3n-\nu}+ \rho_*^{2-n} \, 
\sup_j \|\Phi_j\|_{2, \alpha})
\]
and
\[
\|{\cal M} (F_2) - {\cal N} (F_1)\|_{\cal E}\leq \frac{1}{2}\, \|F_2 -F_1\|_{\cal E}
\]
for all $F_1, F_2 \in {\cal E}_{\rho_*}$ such that $\|F_i\|_{\cal E} \leq c_0 \, (\e^3 \, 
\rho_*^{1-3n- \nu} + \rho_*^{2-n} \, \sup_j \|\Phi_j\|_{2, \alpha})$.

\medskip

In particular ${\cal M}$ has a unique fixed point in this ball.
\end{proposition}
{\bf Proof :} The estimate for ${\cal M} (0)$ follows from the result of Lemma~\ref{le:estimh}
. The other estimate is left to the reader. \hfill $\Box$

\section{Minimal $n$-submanifolds close to the graph of Green's function, parameterized by 
their boundaries}

Again, we give a summary of what we have obtained in the last sections. Assume that $x_1, 
\ldots, x_N \in {\R}^n$,  $\alpha_1, \ldots, \alpha_N \in {\R}$, ${\cal A}_0 \in M_n ({\R})$
and ${\cal R}_1, \ldots, {\cal R}_N \in O(n)$ are fixed. Given $\kappa >0$, $\rho_* >0$,
there exists $\e_0 >0$ such that for all $\alpha_j \in [\frac{1}{\kappa}, \kappa]$, for all
$\e \in (0, \e_0)$ and all $\Phi_j \in {\cal C}^{2, \alpha}(S^{n-1}; {\R}^n)$ which is
orthogonal to ${\cal R}_j \, \Theta$ in the $L^2$ sense and which satisfies
\[
\|\Phi_j\|_{2, \alpha} \leq \kappa \, \rho_*^2 \, \e,
\]
we have obtained a minimal $n$-submanifold which is ${\cal C}^{2, \alpha}_{\nu, \mu}$ 
close to the $n$-plane 
\[
x \longrightarrow x + i  \, \e \, {\cal A}_0 \, x .
\]
and which, up to some translation, can be parameterized, in some neighborhood of each 
of its boundary as the graph of
\[
(r, \theta) \longrightarrow  {\cal R}_ {j_0}  \left( \e \, \beta_{j_0} \, r^{1-n} +  
\frac{\e}{\omega_n}\, \left( \sum_{j\neq j_0} \gamma_{j_0j} \alpha_j - \lambda_{j_0}\right) r
\, \Theta  + J_{j_0}  +  J+ J^\perp \right)
\]
where we have set
\[
 r \, \Theta : = x-x_{j_0}
\]
and where $J_{j_0}$ is the unique harmonic extension of the boundary data $\Phi_j$ 
outside  $B_{\rho_*}$ which tends to $0$ at $\infty$ and where $J_{j_0}$ satisfies
\[
\|J_{j_0}(\rho_*, \cdot)\|_{2, \alpha}+ \|\rho_* \, \del_r J_{j_0} (\rho_*, \dot)\|_{1, 
\alpha} \leq c \, \e  \, \rho_*^2 + c_\kappa \,   \e \, \rho_*^{3-n -\nu}.
\]
and 
\[
\|J^\perp_{j_0} (\rho_*, \cdot)\|_{2, \alpha}+ \|\rho_* \, \del_r J^\perp_{j_0} (\rho_*, 
\dot)\|_{1, \alpha} \leq c \e \, \rho_* .
\]
In addition $J^\perp_{j_0} (r, \cdot)$ is orthogonal to $\Theta$ in the $L^2$ sense on 
$S^{n-1}$ while $J_{j_0} (r, \cdot)$ is collinear to $\Theta$. Observe, and this is
important, that the constants $c>0$ do not depend on $\kappa$ while $c_\kappa >0$ does.

\medskip

We denote by $\Pi_\e  ((\alpha_j)_j, (\Phi_j)_j )$ this submanifold. Though this depends on 
$ ({\cal R}_j)_j $, on $ (x_j)_j$ and on $\rho_*$, we do not write this dependence in the
notation.

\medskip

For further use it will be convenient to define 
\begin{definition}
For all $\Phi \in {\cal C}^{2, \alpha} (S^{n-1}; {\R}^n)$, we define ${\cal P}_{ext}(\Phi)$ 
to be equal to $\del_r F$ where $F$ is the unique harmonic extension of $\Phi$ outside  the
unit ball of ${\R}^n$, which tends to $0$ at $\infty$.
\end{definition}

\section{Gluing procedure}

Assume that $x_1, \ldots, x_{N} \in {\R}^n$, ${\cal A}_0 \in M_n ({\R})$ and ${\cal R}_1, 
\ldots, {\cal R}_{N}\in O(n)$ are fixed in such a way that (H1), (H2) and (H3) hold.  By
assumption, we can  choose  $( \alpha^*_1, \ldots \alpha^*_{N})\in {\R}^N_{+}$ such that 
\[
\sum_j \gamma_{jj'}\, \alpha^*_{j'} = \lambda_j.
\]
We set 
\[
\beta^*_j = \alpha^*_{j}
\]
Finally we fix $\rho_* >0$ small enough and $\kappa >0$ large enough. By the previous 
analysis there exists $\e_0 >0$ such that, for all $\e \in (0, \e_0)$ and for all $\alpha_j,
\beta_j \in [\frac{1}{\kappa,}, \kappa]$, for all $\Phi_j, \tilde \Phi_j \in {\cal C}^{2,
\alpha}(S^{n-1}; {\R}^n)$ which are orthogonal to $\Theta$ in the $L^2$ sense on $S^{n-1}$, 
we can define $\Pi_{\e, ({\cal R}_j)_j , (x_j)_j,  \rho_*} ((\alpha_j)_j, (\Phi_j)_j )$ this
submanifold. This submanifold is, near each of its boundaries, the graph of 
\[
(r, \theta) \longrightarrow  {\cal R}_ {j_0}  \left( \e \, \beta_{j_0} \, r^{1-n} 
\Theta +  \frac{\e}{\omega_n}\, \sum_{j\neq j_0} \gamma_{j_0j} \, (\alpha_j- \alpha_j^*) \, 
r \, \Theta  + J_{j_0}  +  J+ J^\perp \right)
\]
where we have set
\[
 r \, \Theta : = x-x_{j_0}
\]
and where $J_{j_0}$ is the unique harmonic extension of the boundary data $\Phi_j$ outside  
$B_{\rho_*}$ which tends to $0$ at $\infty$ and where both $J$ and $J^\perp$  satisfies
\begin{equation}
\|J_{j_0} (\rho_*, \cdot)\|_{2, \alpha}+ \|\rho_* \, \del_r J_{j_0} (\rho_*, \dot)\|_{1, 
\alpha} \leq c \, \e  \, \rho_*^2 + c_\kappa \,   \e \, \rho_*^{3-n -\nu}.
\label{eq:1}
\end{equation}
and 
\begin{equation}
\|J^\perp_{j_0} (\rho_*, \cdot)\|_{2, \alpha}+ \|\rho_* \, \del_r J^\perp_{j_0} (\rho_*, 
\dot)\|_{1, \alpha} \leq c \, \e \, \rho_* .
\label{eq:2}
\end{equation}
In addition $J^\perp_{j_0} (r, \cdot)$ is orthogonal to $\Theta$ in the $L^2$ sense on 
$S^{n-1}$ while $J_{j_0}(r, \cdot)$ is collinear to $\Theta$. 

\medskip

We can also define a minimal $n$-submanifold $H_{\e, {\cal R}_j} (\beta_j, \tilde \Phi_j)$
 which, once translated by the vector $x_j$ will be denoted by $ H_{\e} (\beta_j, \tilde
\Phi_j)$. Its boundary is parameterized by 
\[
(r, \theta)  \longrightarrow   \e \, \beta  \, r^{1-n} \, \Theta + F_j + \tilde{F}_j,
\]
where $\tilde{F}_j$ depends smoothly on $\beta_j$ and $\Phi_j$ and satisfies
\begin{equation}
\| \tilde F (\rho_*, \cdot) \|_{{\cal C}^{2, \alpha} (S^{n-1})} + \| \rho_* \, \del_r 
\tilde F (\rho_*, \cdot) \|_{{\cal C}^{1, \alpha} (S^{n-1})}\leq c_\kappa  \, \left( \e^3
\rho_*^{1-3n} + \e \, \rho_* \, (\e^{-\frac{1}{n}} \,  \rho_*)^{\frac{2-n}{2}-\delta}\right) ,
\label{eq:3}
\end{equation}
and where $F_j$ is the harmonic extension of $\tilde{\Phi}_j$ in $B_{\rho_*}$..

\medskip

Our aim will now be to find $\Phi_j$, $\tilde \Phi_j$, $\alpha_j$ and $\beta_j$ in such
 a way that 
\[
M_\e : =  \Pi_\e  ((\alpha_j)_j , (\Phi_j)_j ) \cup_{j'=1}^N H_\e (\beta_{j' }, \Phi_{j'}), 
\]
is a ${\cal C}^{1}$ hypersurface. 
 
\medskip

Writing that the boundary of these submanifolds coincide yields the following system of
 equations
\begin{equation}
\begin{array}{rlllll}
\Phi_j - \tilde{\Phi}_j & = &  ( P_0 \, F_j - J_j^\perp ) (\rho_*, \cdot) \\[3mm]
\ds \e \, ( \beta_{j_0}- \alpha_{j_0}) \, \rho_*^{1-n} \Theta +  \frac{\e}{\omega_n}\,  
\sum_{j\neq j_0} \gamma_{j_0j} \, (\alpha_j - \alpha^*_j) \,  \rho_* \, \Theta  &  = &   (
(I-P_0) \, F_j  - J_j ) (\rho_*, \cdot)
\end{array}
\label{eq:c1}
\end{equation}
where the first equation corresponds to the projection over the space of functions 
orthogonal to $\Theta$ in the $L^2$ sense on $S^{n-1}$ and where the second equation
corresponds to the orthogonal projection over the space of functions spanned by $\Theta$.

\medskip

Writing that the conormals at the boundaries coincides yields the following system of
 equations
\begin{equation}
\begin{array}{rlllll}
{\cal P}_{ext} \Phi_j - {\cal P}_{int} \tilde{\Phi}_j & = &  \rho_* \, ( P_0 \, \del_r F_j  -  \del_r J_j^\perp  ) (\rho_*, \cdot) \\[3mm]
\ds  (1-n) \, \e \, ( \beta_{j_0} - \alpha_{j_0} )\, \rho_*^{1-n} +  \frac{\e}{\omega_n}\,
 \sum_{j\neq j_0} \gamma_{j_0j} \, (\alpha_j - \alpha^*_j) \, \rho_* \, \Theta  & = &   
\rho_* \, ( (I - P_0) \, \del_r F_j - \del_r J_j) (\rho_*, \cdot)
\end{array}
\label{eq:c2}
\end{equation}
where the first equation corresponds to the projection over the space of functions 
orthogonal to $\Theta$ in the $L^2$ sense on $S^{n-1}$ and where the second equation
corresponds to the orthogonal projection over the space of functions spanned by $\Theta$.

\medskip
We will now use the well known result
\begin{lemma}
The mapping ${\cal P}_{ext}-{\cal P}_{int}$ is an isomorphism from ${\cal C}^{2, \alpha}
(S^{n-1}; {\R}^n)$ into ${\cal C}^{1, \alpha}(S^{n-1}; {\R}^n)$. 
\end{lemma}
{\bf Proof :}
Observe that both ${\cal P}_{int}$ and ${\cal P}_{ext}$ are self-adjoint first order
 pseudodifferential operator which are elliptic, with principal symbols $|\xi|$ and $-|\xi|$,
respectively, hence the difference is also elliptic and semibounded. This means that ${\cal
P}_{ext}- {\cal P}_{int}$ has discrete spectrum, and thus we need only prove that it is
injective. The invertibility in H\"older spaces then follows by standard regularity theory. 

\medskip

To prove this, we argue by contradiction. Assume that ${\cal P}_{ext}-{\cal P}_{int}$ 
is not injective. Then, there would exist some function $\Phi \in {\cal C}^{2,
\alpha}(S^{n-1}; {\R}^n)$ for which $({\cal P}_{ext}-{\cal P}_{int}) \, \Phi = 0$. We may
extend the Dirichlet data $\Phi$ to a harmonic mapping $F$ on the $B_1$ and also on ${\R}^n
\setminus B_1$. In addition $F$ tends to $0$ at $\infty$. Since $({\cal P}_{ext}-{\cal
P}_{int}) \, \Phi = 0$, $F$ is ${\cal C}^1$ and hence is ${\cal C}^\infty$ and tends to $0$
at $\infty$. Thus $F\equiv 0$. This completes the proof of the injectivity of ${\cal
P}_{ext}-{\cal P}_{int}$. \hfill $\Box$

\medskip
We set 
\[
{\cal F}^\alpha: =[{\cal C}^{2, \alpha}(S^{n-1}; {\R}^n) \times {\R}]^{2N}
\]
endowed with the product norm. Using this result, the previous system of equations 
(\ref{eq:c1}) and (\ref{eq:c2}) reduces to
\[
((\Phi_j)_j, (\tilde \Phi_j)_j, (\e \, (\alpha_j-\alpha_j^*))_j, (\e \, (\beta_j- 
\beta_j^*))_j) = {\bf C}  \left( (\Phi_j)_j, (\tilde \Phi_j)_j, (\e \,
(\alpha_j-\alpha_j^*))_j, (\e \, ( \beta_j -\beta_j^*))_j \right) , 
\]
where the nonlinear mapping ${\bf C}$ satisfies
\[
\| {\bf C} \left( (\Phi_j)_j, (\tilde \Phi_j)_j, (\e \, \alpha_j)_j, (\e \, \beta_j)_j 
\right) \|_{\cal F}\leq c \, \rho_* \, \e
\]
for some constant which does not depend on $\kappa$, provided $\e$ is chosen small 
enough, say $\e \in (0, \e_0)$. This last claim is a simple consequence of
(\ref{eq:1})-(\ref{eq:3}).

\medskip

We denote by ${\cal B}^{\alpha}_{\kappa}$ the ball of radius $\kappa \, \rho_* \, \e$ 
in ${\cal F}^\alpha$. It follows from our previous analysis that, for fixed $\kappa >0$ large
enough, the mapping ${\bf C}$ is well defined in ${\cal B}^\alpha_{\kappa}$ provided the
parameter $\e$ is small enough. 

\medskip

This zero of ${\bf C}$ produces a ${\cal C}^{1, \alpha}$ minimal $n$-submanifold  
$M_\e$. It is then a simple exercise to see, thanks to regularity theory, that $M_\e$  is in
fact a ${\cal C}^\infty$ minimal hypersurface with $N+1$ ends. 

\medskip

To conclude, we want to use Schauder's fixed point Theorem which will ensure the existence 
of at least one fixed point  of ${\bf C}$ in ${\cal B}^\alpha_{\kappa}$. However, since ${\bf
C}$ is not compact it is not possible to apply directly Schauder's Theorem. This is the
reason why we introduce a family of smoothing operators ${\bf D}^{q}$, for all $q >1$, which
satisfy for fixed $0 < \alpha' < \alpha <1$
\[
\| {\bf D}^q \Phi \|_{{\cal C}^{2, \alpha'}(S^{n-1})}\leq c_0 \, \|\Phi \|_{{\cal C}^{2, 
\alpha}(S^{n-1})} \qquad \qquad \| {\bf D}^q \Phi \|_{{\cal C}^{2, \alpha}(S^{n-1})}\leq c_0
\, q^{\alpha -\alpha'} \,  \|\Phi \|_{{\cal C}^{2, \alpha'}(S^{n-1})} ,
\]
and 
\begin{equation}
\| \Phi - {\bf D}^q \Phi \|_{{\cal C}^{2, \alpha'}(S^{n-1})}\leq c_0 \, q^{\alpha' - 
\alpha} \, \| \Phi \|_{{\cal C}^{2, \alpha}(S^{n-1})}.
\label{eq:cde}
\end{equation}
for some constant $c_0>0$ which does not depend on $q>1$. The existence of such smoothing
 operators is available in \cite{Ali-Ger}, Proposition 1.6, page 97.  To keep the notation
short, we use the same notation for the smoothing operator defined on ${\cal F}^\alpha$ and
acting on all function spaces.

\medskip

Now we fix $\kappa > 0$ large enough. For all $q >1$, we may apply Schauder's fixed point
 Theorem to ${\bf D}^q \, {\bf C}$ to obtain the existence of $P_q$ fixed point of ${\bf D}^q
\, {\bf C}$ in ${\cal B}^\alpha_{\kappa}$, provided $\e$ is chosen small enough, say $\e \in
(0, \e_0]$.

\medskip

Since $P_q$ has norm bounded uniformly in $q$, we may extract a sequence $q_j \rightarrow
 +\infty $ such that $P_{q_j}$ converges in ${\cal F}^{\alpha'}$ for some fixed $\alpha' <
\alpha$. Thanks to the continuity of ${\bf C}$  (with respect to the ${\cal C}^{2, \alpha'}$
and ${\cal C}^{1, \alpha'}$ topology) and also to (\ref{eq:cde}), the limit of this sequence
is a fixed point of the mapping ${\bf C}$ and hence, produces a zero of ${\bf C}$, for all
$\e \in (0, \e_0]$. This completes our proof of the Theorem.

\noindent{Claudio Arezzo} \newline
Dipartimento di Matematica \newline
Universit\`{a} di Parma \newline
Via M. D'Azeglio 85 \newline
43100, Parma \newline
Italy

\noindent{e-mail: claudio.arezzo@unipr.it}

\vspace{3mm}

\noindent{Frank Pacard} \newline
Centre de Math\'{e}matiques - Facult\'{e} de Sciences et Technologie \newline
Universit\'{e} Paris XII - Val de Marne \newline
61, Avenue du G\'{e}n\'{e}ral de Gaulle \newline
94 010 Creteil Cedex \newline
France

\noindent{e-mail: pacard@univ-paris12.fr}

\end{document}